\documentclass[english,10pt]{article}
\usepackage[T1]{fontenc}
\usepackage[latin9]{inputenc}
\usepackage{geometry}
\usepackage{textcomp}
\usepackage{amsthm}
\usepackage{amsmath}
\usepackage{amssymb}
\usepackage[all]{xy}

\theoremstyle{plain}
\newtheorem{thm}{Theorem}[subsection]
  \theoremstyle{definition}
  \newtheorem{defn}[thm]{Definition}
  \theoremstyle{plain}
  \newtheorem{conjecture}[thm]{Conjecture}
  \theoremstyle{plain}
  \newtheorem{lem}[thm]{Lemma}
  \theoremstyle{plain}
  \newtheorem{prop}[thm]{Proposition}

\usepackage{mathrsfs}
\usepackage{hyperref}
\usepackage{xy}
\hypersetup{colorlinks=true,citecolor=blue}

\newcommand{\sch}{\mathscr{H}}
\newcommand{\scx}{\mathscr{X}}

\newcommand{\scz}{\mathscr{Z}}
\newcommand{\scw}{\mathscr{W}}
\newcommand{\scm}{\mathscr{M}}
\newcommand{\scv}{\mathscr{V}}
\newcommand{\scc}{\mathscr{C}}
\newcommand{\sco}{\mathscr{O}}
\newcommand{\srr}{\mathscr{R}}
\newcommand{\scd}{\mathscr{D}}
\newcommand{\sca}{\mathscr{A}}

\newcommand{\scl}{\mathscr{L}}

\newcommand{\scg}{\mathscr{G}}

\usepackage[T2A,T1]{fontenc}
\DeclareSymbolFont{cyrillic}{T2A}{cmr}{m}{n}
\DeclareMathSymbol{\Sha}{\mathalpha}{cyrillic}{216}

\usepackage{babel}

\begin{document}

\title{Iwasawa theory of overconvergent modular forms, I: Critical \emph{p}-adic
\emph{L}-functions}

\author{by David Hansen\thanks{Department of Mathematics, Columbia University, 2990 Broadway, New
York NY 10027; hansen@math.columbia.edu}}
\maketitle
\begin{abstract}
We construct an Euler system of \emph{p-}adic zeta elements over the
eigencurve which interpolates Kato's zeta elements over all classical
points. Applying a big regulator map gives rise to a purely algebraic
construction of a two-variable \emph{p-}adic \emph{L-}function over
the eigencurve. As a first application of these ideas, we prove the
equality of the \emph{p-}adic \emph{L-}functions associated with a
critical-slope refinement of a modular form by the works of Bellaïche/Pollack-Stevens
and Kato/Perrin-Riou.

\tableofcontents{}
\end{abstract}

\section{Introduction}

\subsection{Euler systems and zeta elements}

The Birch---Swinnerton-Dyer conjecture and its generalizations relate
special values of \emph{L-}functions with various arithmetically defined
groups. These conjectures are among the most fascinating and difficult
problems in modern number theory. One of the most powerful techniques
for attacking these problems is the theory of \emph{Euler systems,
}families of Galois cohomology classes with very special properties
 \cite{KoEuler,Rubook}. Euler systems forge a strong link between
special values of \emph{L-}functions and Galois cohomology, and every
known Euler system has led to very deep results: among other successes,
we mention Rubin's proof of the Iwasawa main conjectures for imaginary
quadratic fields \cite{RuImquad}, Kolyvagin's results on the Birch---Swinnerton-Dyer
conjecture \cite{Kofinite}, and Kato's work on the Bloch-Kato conjecture
for elliptic modular forms \cite{BK,Kato}.

In this article, we construct Euler systems associated with finite-slope
overconvergent modular forms, and with rigid analytic families of
such forms. Since many classical modular forms are finite-slope overconvergent,
our constructions have consequences for the Iwasawa theory of classical
modular forms and elliptic curves. In particular, Theorem \ref{critical}
below settles the long-standing question of comparing the analytic
and algebraic \emph{p-}adic \emph{L-}functions attached to a critical-slope
refinement of a modular form by Pollack-Stevens and Kato. The proof
of this theorem follows naturally from our construction of a {}``two-variable
algebraic \emph{p-}adic \emph{L-}function\emph{'' }over the eigencurve.
These ideas also lead to some partial results towards a two-variable
main conjecture, which we defer to a second article \cite{Htwovar};
it seems likely that our constructions will have further applications.
In the context of Hida families, many of the results in this paper
have been established by Ochiai in a beautiful series of articles \cite{Och1,Och2,Och3},
which gave us the confidence to attempt this project. We also mention
the related work of Wang \cite{W}, who has proven a lovely two-variable
interpolation of Kato's explicit reciprocity law.

Our point of departure is the monumental work of Kato \cite{Kato}.
Speaking very roughly, we are able to interpolate Kato's constructions
of \emph{p-}adic zeta elements and Euler systems over the Coleman-Mazur-Buzzard
eigencurve. In the course of doing this, we are led to several intermediate
results which may be of interest: in particular, we put canonical
integral structures on Coleman families and on the universal Galois
representations over them.

Before stating our results, we recall the eigencurve and define some
distinguished classes of points on it. Fix an odd prime $p$, an algebraic
closure $\overline{\mathbf{Q}_{p}}$, and an isomorphism $\mathbf{C}\overset{\sim}{\to}\overline{\mathbf{Q}_{p}}$.
Fix an integer $N\geq1$ prime to $p$, and let $\mathbf{T}$ be the
polynomial algebra over $\mathbf{Z}$ generated by the operators $T_{\ell},\ell\nmid Np$,
$U_{p}$ and $\left\langle d\right\rangle ,d\in(\mathbf{Z}/N\mathbf{Z})^{\times}$.
Set $\mathfrak{W}=\mathrm{Spf}(\mathbf{Z}_{p}[[\mathbf{Z}_{p}^{\times}]])$,
and let $\scw=\mathfrak{W}^{\mathrm{rig}}$ be the rigid analytic
space of characters of $\mathbf{Z}_{p}^{\times}$ together with its
universal character $\chi_{\scw}:\mathbf{Z}_{p}^{\times}\to\sco(\scw)^{\times}$;
we embed $\mathbf{Z}$ in $\scw(\mathbf{Q}_{p})$ by mapping $k$
to the character $t\mapsto t^{k-2}$. For any $\lambda\in\scw(\overline{\mathbf{Q}_{p}})$
we (slightly abusively) write \[
M_{\lambda}^{\dagger}(\Gamma_{1}(N))\subset\overline{\mathbf{Q}_{p}}[[q]]\]
for the space of $q$-expansions of overconvergent modular forms of
weight $\lambda$ and tame level $N$. Coleman-Mazur \cite{CMeigencurve}
and Buzzard \cite{BuEigen} constructed a rigid analytic curve $\scc(N)$
equipped with a flat, locally finite morphism $w:\scc(N)\to\scw$
and an algebra homomorphism $\phi:\mathbf{T}\to\sco(\scc(N))$ such
that $\scc(N)(\overline{\mathbf{Q}_{p}})$ parametrizes the set of
overconvergent eigenforms of finite slope and tame level $N$, with
$x\in\scc(N)(\overline{\mathbf{Q}_{p}})$ corresponding to the generalized
eigenspace\begin{eqnarray*}
M_{w(x)}^{\dagger}(\Gamma_{1}(N))_{\ker\phi_{x}} & = & \left\{ f\in M_{w(x)}^{\dagger}(\Gamma_{1}(N))\mid(T-\phi_{x}(T))^{n}f=0\,\forall T\in\mathbf{T}\,\mathrm{and}\, n\gg0\right\} \\
 & \subset & M_{w(x)}^{\dagger}(\Gamma_{1}(N)).\end{eqnarray*}
When this eigenspace is one-dimensional over $\overline{\mathbf{Q}_{p}}$
we write $\mathbf{f}_{x}$ for the canonically normalized generator.
For any $x\in\scc(N)(\overline{\mathbf{Q}_{p}})$, the theory of pseudorepresentations
yields a continuous two-dimensional $G_{\mathbf{Q}}$-representation
$V_{x}$ over the residue field $E_{x}$ of $x$ with tame Artin conductor
dividing $N$ such that for all primes $\ell\nmid Np$ we have\begin{eqnarray*}
\mathrm{tr}\mathrm{Frob}_{\ell}|V_{x} & = & \phi_{x}(T_{\ell}),\\
\det\mathrm{Frob}_{\ell}|V_{x} & = & \ell\chi_{\scw,w(x)}(\ell)\phi_{x}(\left\langle \ell\right\rangle )\end{eqnarray*}
with $\mathrm{Frob}_{\ell}\in G_{\mathbf{Q}}$ a geometric Frobenius.
The function $\alpha=\phi(U_{p})\in\sco(\scc(N))$ is nonvanishing
on $\scc(N)$ and plays a distinguished role in the theory; in particular,
the image of $\scc(N)$ in $\scw\times\mathbf{A}^{1}$ under the map
$x\mapsto(w(x),1/\alpha(x))$ is a hypersurface over which $\scc(N)$
is finite.

Given an integer $M|N$, we say a point $x\in\scc(N)(\overline{\mathbf{Q}_{p}})$
is \emph{crystalline of conductor $M$ }if there is a (necessarily
unique) normalized cuspidal elliptic newform \[
f_{x}=\sum_{n\geq1}a_{n}(f_{x})q^{n}\in S_{k}(\Gamma_{1}(M))\]
of weight $k=k_{x}\geq2$ and nebentype $\varepsilon_{f_{x}}$ such
that $\phi_{x}(T_{\ell})=a_{\ell}(f_{x})$ and $\phi_{x}(\left\langle \ell\right\rangle )=\varepsilon_{f_{x}}(\ell)$
for all $\ell\nmid Np$ and furthermore $\alpha_{x}$ is a root of
the Hecke polynomial $X^{2}-a_{p}(f_{x})+p^{k_{x}-1}\varepsilon_{f_{x}}(p)$.\emph{}%
\footnote{These are precisely the points for which the Galois representation
$V_{x}$ has tame Artin conductor $M$ and is crystalline with distinct
nonnegative Hodge-Tate weights at $p$.%
}It's not hard to see that any newform of weight $\geq2$ and level
dividing $N$ occurs this way. Given any point $x\in\scc(N)(\overline{\mathbf{Q}_{p}})$
with $w(x)=k,k\geq2$, we say $x$ is \emph{noncritical }if the inclusion
of generalized eigenspaces\[
M_{k}(\Gamma_{1}(N)\cap\Gamma_{0}(p))_{\ker\phi_{x}}\subset M_{k}^{\dagger}(\Gamma_{1}(N))_{\ker\phi_{x}}\]
is an equality, and \emph{critical }otherwise.
\begin{defn}
\emph{A point $x\in\scc(N)(\overline{\mathbf{Q}_{p}})$ is }\textbf{noble}
\emph{if it is crystalline of conductor $N$, and furthermore $x$
is noncritical and $\alpha_{x}$ is a simple root of the pth Hecke
polynomial of $f_{x}$.}
\end{defn}
Note that a crystalline point of conductor $N$ is noble exactly when
the generalized eigenspaces in the definition of noncriticality are
honest eigenspaces of common dimension one, in which case the \emph{p-}stabilized
form\[
\mathbf{f}_{x}=f_{x}(q)-p^{k_{x}-1}\varepsilon_{f_{x}}(p)\alpha_{x}^{-1}f_{x}(q^{p})\]
defines the canonical generator of both spaces. Our motivation for
singling out these points is geometric: the eigencurve is smooth and
étale over $\scw$ locally at any noble point \cite{BellCrit}, and
certain sheaves of interest are locally free of minimal rank around
such points. It is conjectured that the final condition in the definition
of nobility always holds \cite{CE}, and that a crystalline point
of conductor $N$ fails to be noble if and only if $f_{x}$ is a CM
form and $v_{p}(\alpha_{x})=k_{x}-1$; in particular, the map\begin{eqnarray*}
\left\{ \mathrm{noble\, points}\right\}  & \to & \left\{ \mathrm{newforms}\,\mathrm{of}\,\mathrm{level}\, N\,\mathrm{and}\,\mathrm{weight}\,\geq2\right\} \\
x & \mapsto & f_{x}\end{eqnarray*}
is conjecturally surjective, one-to-one on the set of CM forms with
$p$ split in the CM field, and two-to-one on all other forms. These
results are all known for points of weight $2$, and in particular,
any point such that $f_{x}$ is associated with a non-CM (modular)
elliptic curve over $\mathbf{Q}$ is noble.

Let $\scc_{0}^{M-\mathrm{new}}(N)$ denote the Zariski closure in
$\scc(N)$ of the crystalline points of conductor $M$; this is a
union of irreducible components of $\scc(N)$, and the complement
of $\cup_{M|N}\scc_{0}^{M-\mathrm{new}}(N)$ in $\scc(N)$ is an open
dense subspace of the union of Eisenstein components of $\scc(N)$.
\begin{defn}
\emph{The }smooth cuspidal eigencurve of level $N$, \emph{denoted
by $\scc_{N}$ or simply by $\scc$},\emph{ is the normalization of
}$\scc_{0}^{N-\mathrm{new}}(N)$.
\end{defn}
The curve $\scc$ is a disjoint union of connected smooth reduced
rigid analytic curves, with a natural morphism $i:\scc\to\scc(N)$.
We will slightly abusively write $w=w\circ i$, $\phi=i^{\ast}\phi$,
$\alpha=i^{\ast}\alpha$ for the natural maps and functions inherited
from those on $\scc(N)$. Any point in the smooth locus of $\scc_{0}^{N-\mathrm{new}}(N)$
determines a unique point in $\scc$, and $i$ is an isomorphism locally
around any such point; in particular, this holds for any noble point\emph{.}
We find it very convenient to introduce the involution $x\mapsto x^{c}$
of $\scc$ given on geometric points by\begin{eqnarray*}
\phi_{x^{c}}(T_{\ell}) & = & \phi_{x}(T_{\ell}\left\langle \ell\right\rangle ^{-1}),\\
\phi_{x^{c}}(U_{p}) & = & \phi_{x}(U_{p}\left\langle p\right\rangle ^{-1}),\\
\phi_{x^{c}}(\left\langle d\right\rangle ) & = & \phi_{x}(\left\langle d^{-1}\right\rangle ).\end{eqnarray*}
Note that on crystalline points we have $f_{x^{c}}=f_{x}^{c}$ where
$f^{c}=f\otimes\varepsilon_{f}^{-1}$ is the complex conjugate of
$f$, and likewise for $\alpha_{x}$.

As we construct it here, the curve $\scc$ carries a rank two locally
free sheaf $\scv$ with a continuous $G_{\mathbf{Q}}$-action unramified
outside $Np$, such that for every noble point $x\in\scc(\overline{\mathbf{Q}_{p}})$
there is an isomorphism $\scv_{x}\simeq V_{f_{x}}(k_{x})\simeq V_{f_{x}^{c}}^{\ast}(1)$
which is moreover realized by a \emph{canonical }isomorphism\[
\scv_{x}\cong\left(H_{\mathrm{\acute{e}t}}^{1}\left(Y_{1}(Np)_{/\overline{\mathbf{Q}}},(\mathrm{sym}^{k_{x}-2}T_{p}E)(2)\right)\otimes_{\mathbf{Z}_{p}}E_{x}\right)[\ker\phi_{x}]\]
of Galois modules. Here $T_{p}E$ is the étale sheaf defined by the
Tate module of the universal elliptic curve over $Y_{1}(Np)$. The
sheaf $\scv$ and its fibers provide canonical étale-cohomological
realizations of the Galois representations associated with the overconvergent
eigenforms parametrized by $\scc$. Following ideas of Pottharst,
we define coherent Galois cohomology sheaves $\sch^{1}(\mathbf{Q}(\zeta_{m}),\scv(-r))$
over $\scc$ whose sections over an admissible affinoid open $U$
are given by\[
\sch^{1}(\mathbf{Q}(\zeta_{m}),\scv(-r))(U)=H^{1}(G_{\mathbf{Q}(\zeta_{m}),Nmp\infty},\scv(U)(-r)).\]
We denote the global sections of this sheaf by $H^{1}(\mathbf{Q}(\zeta_{m}),\scv(-r))$.

As a sample of our results, we state the following theorem.
\begin{thm}
\label{euler1}Notation as above, let $\nu$ be any primitive Dirichlet
character of even order and conductor $A$ prime to $Np$, with coefficient
field $\mathbf{Q}_{p}(\nu)$, and let $r$ be an arbitrary integer.

\emph{i.} There is a collection of global cohomology classes\[
\mathfrak{z}(r,\nu)=\left(\mathfrak{z}_{m}(r,\nu)\in H^{1}(\mathbf{Q}(\zeta_{m}),\scv(-r))\otimes_{\mathbf{Q}_{p}}\mathbf{Q}_{p}(\nu)\right){}_{m\geq1,\,(m,A)=1}\]
such that for any prime $\ell\nmid NA$, we have the corestriction
relation\emph{\[
\mathrm{Cor}_{\mathbf{Q}(\zeta_{m\ell})/\mathbf{Q}(\zeta_{m})}\left(\mathfrak{z}_{\ell m}(r,\nu)\right)=\begin{cases}
\mathfrak{z}_{m}(r,\nu) & \mathrm{if}\,\ell|mp\\
P_{\ell}(\ell^{-r}\sigma_{\ell}^{-1})\cdot\mathfrak{z}_{m}(r,\nu) & \mathrm{if}\ell\nmid mp,\end{cases}\]
}where $\sigma_{\ell}\in\mathrm{Gal}(\mathbf{Q}(\zeta_{m})/\mathbf{Q})$
denotes an arithmetic Frobenius at $\ell$ and $P_{\ell}(X)\in\sco(\scc)[X]$
is the polynomial\begin{eqnarray*}
P_{\ell}(X) & = & 1-\phi(T_{\ell}\left\langle \ell\right\rangle ^{-1})X+\phi(\left\langle \ell\right\rangle ^{-1})\chi_{\scw}(\ell)\ell X^{2}\\
 & = & \mathrm{det}\left(1-X\cdot\mathrm{Frob}_{\ell}\right)|\scv^{\ast}(1).\end{eqnarray*}

\emph{ii.} The class $\mathfrak{z}_{m}(r,\nu)$ is unramified outside
the primes dividing $p$.

\emph{iii.} If $x$ is a noble point such that either $k_{x}>2$ or
$k_{x}=2$ and $L(1,f_{x}\otimes\nu)\neq0$, the specialization\[
\mathfrak{z}(r,\nu)_{x}=\left(\mathfrak{z}_{m}(r,\nu)_{x}\in H^{1}(\mathbf{Q}(\zeta_{m}),\scv_{x}(-r))\otimes_{\mathbf{Q}_{p}}\mathbf{Q}_{p}(\nu)\right)_{m\geq1,\,(m,A)=1}\]
is nonzero.
\end{thm}
The corestriction relation in part i. here is the famous Euler system
relation.

We also construct classes in Iwasawa cohomology which interpolate
the classes in Theorem \ref{euler1} for varying $r$ and $m$. More
precisely, set $\Gamma_{m}=\mathrm{Gal}(\mathbf{Q}(\zeta_{mp^{\infty}})/\mathbf{Q}(\zeta_{m}))$
and let $\Lambda_{m}$ be the sheaf of rings over $\scc$ defined
by \[
\Lambda_{m}(U)=\left(\sco(U)^{\circ}\widehat{\otimes}\mathbf{Z}_{p}[[\Gamma_{m}]]\right)[\frac{1}{p}]\]
for $U$ an affinoid open (we notate these and other allied objects
in the case $m=1$ by dropping $m$).%
\footnote{If $A$ is a reduced affinoid algebra and $R$ is an adic Noetherian
$\mathbf{Z}_{p}$-algebra which is \emph{p-}adically complete, then
$A^{\circ}\widehat{\otimes}R$ denotes the $1\otimes I$-adic completion
of $A^{\circ}\otimes_{\mathbf{Z}_{p}}R$ where $I$ is the largest
ideal of definition of $R$ (so in particular $p\in I$). %
} We construct a sheaf $\sch_{\mathrm{Iw}}^{1}(\mathbf{Q}(\zeta_{m}),\scv(-r))$
of $\Lambda_{m}$-modules characterized by the equality\[
\sch_{\mathrm{Iw}}^{1}(\mathbf{Q}(\zeta_{m}),\scv(-r))(U)=H^{1}(G_{\mathbf{Q}(\zeta_{m}),Nmp\infty},\scv(U)(-r)\otimes_{\sco(U)}\Lambda_{m}^{\iota}(U)).\]
For any integer $r$ and any integer $j\geq0$, there is a natural
morphism of sheaves\[
\theta_{r,j}:\sch_{\mathrm{Iw}}^{1}(\mathbf{Q}(\zeta_{m}),\scv(-r))\to\sch^{1}(\mathbf{Q}(\zeta_{mp^{j}}),\scv(-r))\]
induced by the surjections $\Lambda_{m}(U)\twoheadrightarrow\sco(U)\otimes_{\mathbf{Q}_{p}}\mathbf{Q}_{p}[\mathrm{Gal}(\mathbf{Q}(\zeta_{mp^{j}})/\mathbf{Q}(\zeta_{m}))]$
together with Shapiro's lemma. Let $H_{\mathrm{Iw}}^{1}$ denote the
global sections of $\sch_{\mathrm{Iw}}^{1}$.
\begin{thm}
\label{iwasawa1}There are naturally defined cohomology classes\[
\mathfrak{Z}_{m}(r,\nu)\in H_{\mathrm{Iw}}^{1}(\mathbf{Q}(\zeta_{m}),\scv(-r))\otimes_{\mathbf{Q}_{p}}\mathbf{Q}_{p}(\nu)\]
whose image under $\theta_{r,j}$ equals $\mathfrak{z}_{mp^{j}}(r,\nu)$
for all $r\in\mathbf{Z}$ and all $j\in\mathbf{Z}_{\geq1}$, and such
that under the natural twisting isomorphisms\[
\mathrm{Tw}_{r_{1}-r_{2}}:\sch_{\mathrm{Iw}}^{1}(\mathbf{Q}(\zeta_{m}),\scv(-r_{1}))\overset{\sim}{\to}\sch_{\mathrm{Iw}}^{1}(\mathbf{Q}(\zeta_{m}),\scv(-r_{2}))\]
we have $\mathrm{Tw}_{r_{1}-r_{2}}(\mathfrak{Z}_{m}(r_{1},\nu))=\mathfrak{Z}_{m}(r_{2},\nu)$.
\end{thm}
Theorem \ref{euler1}, while aesthetically pleasing, is somewhat useless
in applications, since the Euler system machinery deals in integral
structures. In fact Theorems \ref{euler1} and \ref{iwasawa1} follow
from our main technical result, which takes integral structures into
account and which we now describe. Given $k\geq2$ and $M\geq1$,
set\[
V_{k}^{\circ}(M)=H_{\acute{\mathrm{e}}\mathrm{t}}^{1}\left(Y_{1}(M)_{/\overline{\mathbf{Q}}},\mathrm{sym}^{k-2}T_{p}E\right)(2-k).\]
This space (denoted $V_{k,\mathbf{Z}_{p}}(Y_{1}(M))$ in \cite{Kato})
admits commuting actions of $G_{\mathbf{Q}}$ and a suitable algebra
of Hecke operators, and its eigenspace for the Hecke eigenvalues of
a newform $f$ of weight $k$ and level $M$ provides a canonical
realization of the Galois representation associated with $f$. Generalizing
ideas of Stevens and others, given a suitable {}``small'' formal
subscheme $\mathfrak{U}\subset\mathfrak{W}$ and a sufficiently large
integer $s$, we define a module $\mathcal{D}_{\mathfrak{U}}^{s,\circ}$
of $\sco(\mathfrak{U})$-valued locally analytic distributions on
$\mathbf{Z}_{p}$ with a continuous action of $\Gamma_{0}(p)$. Here
{}``small'' means basically that $\mathfrak{U}=\mathrm{Spf}(R)$
where $R=\sco(\mathfrak{U})$ is a $\mathbf{Z}_{p}$-flat and module-finite
$\mathbf{Z}_{p}[[X_{1},\dots,X_{d}]]$-algebra - in particular, the
Berthelot generic fiber $\mathfrak{U}^{\mathrm{rig}}\subset\scw$
is typically not quasicompact. On the other hand, working with small
opens like this allows us to construct a decreasing filtration on
the module $\mathcal{D}_{\mathfrak{U}}^{s,\circ}$ by sub-$\sco(\mathfrak{U})$-modules
\[
\mathcal{D}_{\mathfrak{U}}^{s,\circ}=\mathrm{Fil}^{0}\mathcal{D}_{\mathfrak{U}}^{s,\circ}\supset\mathrm{Fil}^{1}\mathcal{D}_{\mathfrak{U}}^{s,\circ}\supset\cdots\supset\mathrm{Fil}^{i}\mathcal{D}_{\mathfrak{U}}^{s,\circ}\supset\cdots\]
with several extremely felicitous properties: this filtration is $\Gamma_{0}(p)$-stable,
$\mathcal{D}_{\mathfrak{U}}^{s,\circ}/\mathrm{Fil}^{i}\mathcal{D}_{\mathfrak{U}}^{s,\circ}$
is a finite abelian group of exponent $p^{i}$, $\Gamma(p^{s+i})$
acts trivially on $\mathcal{D}_{\mathfrak{U}}^{s,\circ}/\mathrm{Fil}^{i}\mathcal{D}_{\mathfrak{U}}^{s,\circ}$,
and $\mathcal{D}_{\mathfrak{U}}^{s,\circ}$ is separated and complete
for the topology defined by $\mathrm{Fil}^{\bullet}\mathcal{D}_{\mathfrak{U}}^{s,\circ}$.
With these facts in hand we are able to show that the module \[
\mathbf{V}_{\mathfrak{U}}^{s,\circ}(N):=\lim_{\infty\leftarrow i}H_{\acute{\mathrm{e}}\mathrm{t}}^{1}\left(Y_{1}(pN)_{/\overline{\mathbf{Q}}},\mathcal{D}_{\mathfrak{U}}^{s,\circ}/\mathrm{Fil}^{i}\mathcal{D}_{\mathfrak{U}}^{s,\circ}\right)(2)\]
is canonically isomorphic as a Hecke module to the cohomology of the
local system induced by $\mathcal{D}_{\mathfrak{U}}^{s,\circ}$ on
the analytic space $Y_{1}(Np)(\mathbf{C})$. $ $The modules $\mathbf{V}_{\mathfrak{U}}^{s,\circ}(N)$
and variants thereof are the main technical objects in this article.
By construction $\mathbf{V}_{\mathfrak{U}}^{s,\circ}(N)$ admits commuting
Hecke and Galois actions, and we will see that in addition the Galois
action is $p$-adically continuous and unramified outside the primes
dividing $Np$. For any $k\in\mathfrak{U}\cap\mathbf{Z}_{\geq2}$,
there is a natural Hecke- and Galois-equivariant {}``integration''
map\[
i_{k}:\mathbf{V}_{\mathfrak{U}}^{s,\circ}(N)\to V_{k}^{\circ}(Np)(k).\]

\begin{thm}
\textbf{\emph{(cf. Proposition 3.2.1) }}\label{eulermain}For any
integer $r$, any integer $A\geq1$ prime to $p$, any residue class
$a(A)$, any integer $m\geq1$, and any integers $c,d$ with $(cd,6Apm)=(d,Np)=1$,
there is a canonically defined cohomology class\[
_{c,d}\mathbf{\mathfrak{z}}_{N,m}^{s}(\mathfrak{U},r,a(A))\in H^{1}(G_{\mathbf{Q}(\zeta_{m}),Nmp\infty},\mathbf{V}_{\mathfrak{U}}^{s,\circ}(N)(-r))\]
whose specialization under the integration map $i_{k}$ for any $k\in\mathbf{Z}_{\geq2}\cap\mathfrak{U}$
satisfies\[
i_{k}({}_{c,d}\mathbf{\mathfrak{z}}_{N,m}^{s}(\mathfrak{U},r,a(A)))=\phantom{}_{c,d}\mathbf{z}_{1,Np,m}^{(p)}(k,r,k-1,a(A),\mathrm{prime}(mpA))\in H^{1}(G_{\mathbf{Q}(\zeta_{m}),Nmp\infty},V_{k}^{\circ}(Np)(k-r)),\]
where $_{c,d}\mathbf{z}_{1,Np,m}^{(p)}(k,r,r',a(A),S)$ is the p-adic
zeta element defined in §8.9 of\emph{ \cite{Kato}}.
\end{thm}
These classes satisfy the Euler system relations as $m$ varies and
are compatible under changing $\mathfrak{U}$ and $s$, among other
properties. With this result in hand, we prove Theorems \ref{euler1}
and \ref{iwasawa1} by building the eigencurve and all associated
structures from the finite-slope direct summands of the modules $\mathbf{V}_{\mathfrak{U}}^{s}(N)=\mathbf{V}_{\mathfrak{U}}^{s,\circ}(N)[\frac{1}{p}]$.
Precisely, given a pair $(\mathfrak{U},h)$ with $\mathfrak{U}\subset\mathfrak{W}$
as before and $h\in\mathbf{Q}_{\geq0}$, we say this pair is a \emph{slope
datum }if the module $\mathbf{V}_{\mathfrak{U}}^{s}(N)$ admits a
slope-$\leq h$ direct summand $\mathbf{V}_{\mathfrak{U},h}^{s}(N)$
for the action of the $U_{p}$-operator. It turns out that any such
direct summand is module-finite over $\sco(\mathfrak{U})[\frac{1}{p}]$
and preserved by the Hecke and Galois actions. Furthermore, the natural
map $\mathbf{V}_{\mathfrak{U}}^{s+1,\circ}(N)\to\mathbf{V}_{\mathfrak{U}}^{s,\circ}(N)$
induces Hecke- and Galois-equivariant \emph{isomorphisms $\mathbf{V}_{\mathfrak{U},h}^{s+1}(N)\cong\mathbf{V}_{\mathfrak{U},h}^{s}(N)$};
we set $\scv_{\mathfrak{U},h}(N)=\lim_{\leftarrow s}\mathbf{V}_{\mathfrak{U},h}^{s}(N)$.
$ $We construct the full eigencurve $\scc(N)$ by gluing the quasi-Stein
rigid spaces $\scc_{\mathfrak{U}^{\mathrm{rig}},h}(N)$ associated
with the finite $\sco(\mathfrak{U})[\frac{1}{p}]$-algebras\[
\mathbf{T}_{\mathfrak{U},h}(N)=\mathrm{image\,}\mathrm{of}\,\mathbf{T}\otimes_{\mathbf{Z}}\sco(\mathfrak{U})[\tfrac{1}{p}]\,\mathrm{in}\,\mathrm{End}_{\sco(\mathfrak{U})[\frac{1}{p}]}\left(\scv_{\mathfrak{U},h}(N)\right)\]
(or rather, associated with their {}``generic fibers'': we have
$\mathbf{T}_{\mathfrak{U},h}(N)\cong\sco(\scc_{\mathfrak{U}^{\mathrm{rig}},h}(N))^{b}$).
The $\mathbf{T}_{\mathfrak{U},h}(N)$-modules $\scv_{\mathfrak{U},h}(N)$
and $H^{1}(G_{\mathbf{Q}(\zeta_{m}),Nmp\infty},\scv_{\mathfrak{U},h}(N)(-r))$
glue together over this covering into coherent sheaves $\scv(N)$
and $\sch^{1}(\mathbf{Q}(\zeta_{m}),\scv(N)(-r))$, respectively,
and the classes described in Theorem \ref{eulermain} glue into global
sections $\phantom{}_{c,d}\mathfrak{z}_{N,m}(r,a(A))$ of the latter
sheaf. Choosing $\nu$ as above and setting \[
\mu_{c,d}(r,\nu)=(c^{2}-c^{r}\chi_{\scw}(c)^{-1}\nu(c)^{-1}\sigma_{c})(d^{2}-d^{r}\nu(d)\sigma_{d})\in\left(\sco(\scw)\otimes_{\mathbf{Q}_{p}}\mathbf{Q}_{p}(\nu)\right)[\mathrm{Gal}(\mathbf{Q}(\zeta_{m})/\mathbf{Q}],\]
we are able to show that for a judicious choice of $c$ and $d$,
the quantity $\mu_{c,d}$ is a unit in $$\left(\sco(\scw)\otimes_{\mathbf{Q}_{p}}\mathbf{Q}_{p}(\nu)\right)[\mathrm{Gal}(\mathbf{Q}(\zeta_{m})/\mathbf{Q}],$$
and\[
\mathfrak{z}_{N,m}(r,\nu)=\mu_{c,d}(r,\nu)^{-1}\sum_{a\in(\mathbf{Z}/A\mathbf{Z})^{\times}}\nu(a)\phantom{}_{c,d}\mathfrak{z}_{N,m}(r,a(A))\]
is well-defined independently of choosing $c,d$. Finally, we pull
everything back under $i:\scc\to\scc(N)$. The proof of Theorem \ref{iwasawa1}
is only slightly different, although it requires an additional intermediate
result which seems interesting in its own right: we show that the
module\[
\mathbf{V}_{\mathfrak{U},h}^{s,\circ}(N)=\mathrm{image\, of}\,\mathbf{V}_{\mathfrak{U}}^{s,\circ}(N)\,\mathrm{in}\,\mathbf{V}_{\mathfrak{U},h}^{s}(N)\]
is a \emph{finitely presented }$\sco(\mathfrak{U})$-module. Since
$\mathbf{V}_{\mathfrak{U},h}^{s,\circ}(N)[\frac{1}{p}]\cong\mathbf{V}_{\mathfrak{U},h}^{s}(N)$,
this yields natural integral structures on Coleman families. We refer
the reader to §2 and §3 for the details of all these constructions.

\subsection{Critical slope \emph{p}-adic \emph{L}-functions}

We now explain an application of these results to \emph{p-}adic \emph{L-}functions.
To do this, we need a little more notation. For the remainder of the
paper we set\[
\Gamma=\Gamma_{1}=\mathrm{Gal}(\mathbf{Q}_{p}(\zeta_{p^{\infty}})/\mathbf{Q}_{p})\cong\mathrm{Gal}(\mathbf{Q}(\zeta_{p^{\infty}})/\mathbf{Q})\cong\mathbf{Z}_{p}^{\times}.\]
Let $\scx$ be the rigid generic fiber of the formal scheme $\mathrm{Spf}(\mathbf{Z}_{p}[[\Gamma]])$,
so there is a natural bijection of sets (or even of groups)\[
\scx(\overline{\mathbf{Q}_{p}})=\mathrm{Hom}_{\mathrm{cts}}(\Gamma,\overline{\mathbf{Q}_{p}}^{\times}).\]
Given a point $y\in\scx(\overline{\mathbf{Q}_{p}})$ we write $\chi_{y}$
for the corresponding character, and likewise given a character $\chi$
we write $y_{\chi}$ for the corresponding point. Let $\scg$ be the
sheaf of rings over $\scc$ defined by $\scg(U)=\sco(U\times\scx)$,
so $\scg$ is a sheaf of $\Lambda$-algebras. We regard an element
$g\in\scg(U)$ as an analytic function of the two variables $x\in U(\overline{\mathbf{Q}_{p}})$
and $\chi\in\mathrm{Hom}(\Gamma,\overline{\mathbf{Q}_{p}}^{\times})$,
writing $g_{x}$ for its specialization to an element of $\sco(\scx)\otimes_{\mathbf{Q}_{p}}E_{x}$
and $g_{x}(\chi)$ or $g(x,\chi)\in\overline{\mathbf{Q}_{p}}$ for
its value at the point $(x,y_{\chi})$.

Let $f=\sum_{n\geq1}a_{n}(f)q^{n}\in S_{k}(\Gamma_{1}(N))$ be a cuspidal
newform with character $\varepsilon$ and coefficient field $\mathbf{Q}(f)\subset\mathbf{C}$.
Given any element $\gamma\in H^{1}(Y_{1}(N)(\mathbf{C}),\mathrm{sym}^{k-2}(\mathbf{Q}(f)^{2}))$
in the $f$-eigenspace whose projections $\gamma^{\pm}$ to the $\pm1$-eigenspaces
for complex conjugation are nonzero, we define the periods of $f$
by \[
\omega_{f}=\Omega_{f,\gamma}^{+}\gamma^{+}+\Omega_{f,\gamma}^{-}\gamma^{-}\in H_{\mathrm{dR}}^{1}(Y_{1}(N),\mathrm{sym}^{k-2}(\mathbf{C}^{2})).\]
By a fundamental theorem of Eichler and Shimura, for any integer $0\leq j\leq k-2$
and any Dirichlet character $\eta$, the ratio\[
L^{\mathrm{alg}}(j+1,f\otimes\eta)=\frac{j!L(j+1,f\otimes\eta)}{(2\pi i)^{j+1}\tau(\eta)\Omega_{f,\gamma}^{\pm}}\]
is \emph{algebraic, }and in fact lies in the extension $\mathbf{Q}(f,\eta)$;
here the sign here is determined by $(-1)^{j}=\pm\eta(-1)$. 

A \emph{p-}adic \emph{L-}function associated with $f$ interpolates
the special values $L^{\mathrm{alg}}(j+1,f\otimes\eta)$ as $\eta$
varies over Dirichlet characters of $p$-power conductor. More precisely,
let $\scx^{(k)}\subset\scx(\overline{\mathbf{Q}_{p}})$ denote the
set of characters of the form $t\mapsto t^{j}\eta(t)$ with $0\leq j\leq k-2$
and $\eta$ nontrivial of finite order\emph{.} Let $\mathbf{f}$ be
either \emph{p-}stabilization of $f$, with $U_{p}\mathbf{f}=\alpha_{\mathbf{f}}\mathbf{f}$,
and let $\mathbf{Q}_{p}(\mathbf{f})$ be the extension of $\mathbf{Q}_{p}$
generated by $\mathbf{Q}(f)$ and $\alpha_{\mathbf{f}}$. There are
then two natural \emph{p-}adic \emph{L-}functions \[
L_{p,\bullet}(\mathbf{f})\in\mathscr{O}(\scx)\otimes_{\mathbf{Q}_{p}}\mathbf{Q}_{p}(\mathbf{f}),\;\bullet\in\{\mathrm{an},\mathrm{alg}\}\]
associated with $\mathbf{f}$, which we call the \emph{analytic }and
\emph{algebraic p-}adic \emph{L}-functions. Both \emph{L-}functions
satisfy the interpolation formula\[
L_{p,\bullet}(\mathbf{f})(x^{j}\eta(x))=\alpha_{\mathbf{f}}^{-n}p^{(j+1)n}L^{\mathrm{alg}}(j+1,f\otimes\eta^{-1})\]
for all $t^{j}\eta(t)\in\scx^{(k)}$ (here $p^{n}$ is the conductor
of $\eta$) and they both have growth {}``of order at most $v_{p}(\alpha_{\mathbf{f}})$''
in a certain technical sense. After the stimulating initial work of
Mazur---Swinnerton-Dyer \cite{MazurSW}, the analytic \emph{p-}adic
\emph{L-}function was constructed by many people using the theory
of modular symbols in various guises; the cleanest approach is via
Stevens's theory of overconvergent modular symbols, as developed by
Stevens, Pollack-Stevens, and Bellaïche \cite{Strigidsymbs,PSovercvgt,PScrit,BellCrit}.
The algebraic \emph{p-}adic \emph{L-}function was constructed by Kato
as the image of a certain globally defined zeta element in Iwasawa
cohomology under the Perrin-Riou regulator map. 

According to a theorem of Višik \cite{Visik}, any element of $\sco(\scx)$
with growth of order less than $k-1$ is determined uniquely by its
values at the points in $\scx^{(k)}$. When $v_{p}(\alpha_{\mathbf{f}})<k-1$,
this immediately implies the equality $L_{p,\mathrm{an}}(\mathbf{f})=L_{p,\mathrm{alg}}(\mathbf{f})$.
On the other hand, suppose we are in the \emph{critical slope case
}where $v_{p}(\alpha_{\mathbf{f}})=k-1$. Even though the two \emph{p-}adic
\emph{L-}functions agree at all points in $\scx^{(k)}$, Višik's theorem
does not apply, and Pollack and Stevens have pointed out (Remark 9.7
of  \cite{PScrit}) that comparing $L_{p,\mathrm{an}}$ and $L_{p,\mathrm{alg}}$
is a genuine question in this situation. 
\begin{thm}
\label{critical}If $v_{p}(\alpha_{\mathbf{f}})=k-1$ and $V_{f}|G_{\mathbf{Q}_{p}}$
is indecomposable, then $L_{p,\mathrm{an}}(\mathbf{f})=L_{p,\mathrm{alg}}(\mathbf{f})$.
\end{thm}
To the best of our knowledge, this theorem is the first comparison
of two \emph{p-}adic \emph{L-}functions which isn't an immediate consequence
of Weierstrass preparation, Višik's theorem, or some other {}``soft''
result in \emph{p-}adic analysis. When $V_{f}|G_{\mathbf{Q}_{p}}$
is split, $f$ is conjecturally\emph{ }a CM form, in which case this
result is due to Lei-Loeffler-Zerbes \cite{LLZCrit}. Our method is
completely different, and it would be interesting to know if we could
treat the split case without the CM hypothesis.

\emph{Ideally, }the proof of this theorem would proceed as follows.
Let $U$ be a suitably small affinoid neighborhood of the point $x_{\mathbf{f}}$
corresponding to $\mathbf{f}$. The indecomposability of $V_{f}|G_{\mathbf{Q}_{p}}$
implies $x_{\mathbf{f}}$ is a noble point  \cite{BellCrit}. By the
work of many authors (cf. ibid.), any choice of a section $\gamma\in\scv(U)$
as above gives rise to a canonically defined two-variable interpolation
of $L_{p,\mathrm{an}}$, i.e. a function $\mathbf{L}_{p,\mathrm{an}}\in\sco(U\times\scx)$
with the following properties:
\begin{itemize}
\item $\mathbf{L}_{p,\mathrm{an},x}=L_{p,\mathrm{an}}(\mathbf{f}_{x})$
for all noble points $x\in U(\overline{\mathbf{Q}_{p}})$, and
\item $\mathbf{L}_{p,\mathrm{an},x}$ has growth of order at most $v_{p}(\alpha_{x})$
for any\emph{ }$x\in U(\overline{\mathbf{Q}_{p}}).$
\end{itemize}
\emph{Suppose} we could construct an analogous function $\mathbf{L}_{p,\mathrm{alg}}$
interpolating $L_{p,\mathrm{alg}}$ and with exactly the same growth
properties as $\mathbf{L}_{p,\mathrm{an}}$. By the remarks immediately
preceding Theorem \ref{critical}, for any noble point $x$ of noncritical
slope the difference \[
\mathbf{L}_{p,\mathrm{an},x}-\mathbf{L}_{p,\mathrm{alg},x}=L_{p,\mathrm{an}}(\mathbf{f}_{x})-L_{p,\mathrm{alg}}(\mathbf{f}_{x})\]
is identically zero in $\sco(\scx)\otimes_{\mathbf{Q}_{p}}E_{x}$,
so the Zariski-density of such points in $U$ implies the equality
$\mathbf{L}_{p,\mathrm{an}}=\mathbf{L}_{p,\mathrm{alg}}$. Specializing
this equality at our original point $x_{\mathbf{f}}$ implies Theorem
\ref{critical}. Of course, the existence of $\mathbf{L}_{p,\mathrm{alg}}$
is the really nontrivial ingredient in this argument. We don't actually
construct $\mathbf{L}_{p,\mathrm{alg}}$ on the nose in this paper,
but we do construct a function close enough in behavior that a slightly
modified version of the above argument goes through.

To give this construction, we introduce a regulator map which interpolates
the Perrin-Riou regulator map over $\scc$. Let $\srr$ be the sheaf
of relative Robba rings over $\scc$, so $\srr(U)\cong\srr_{\mathbf{Q}_{p}}\widehat{\otimes}\sco(U)$
for any admissible open affinoid $U$. The sheaf $\srr$ carries commuting
$\sco_{\scc}$-linear actions of $\Gamma$ and an operator $\varphi$.
Set $\scv_{0}=\scv(-1)$ and $\scv_{0}^{\ast}=\sch\mathrm{om}_{\sco_{\scc}}(\scv_{0},\sco_{\scc})$,
so $\scv_{0,x}^{\ast}\simeq V_{f_{x}^{c}}$ at any point. Let $\scd^{\ast}=\mathbf{D}_{\mathrm{rig}}^{\dagger}(\scv_{0}^{\ast}|G_{\mathbf{Q}_{p}})$
be the relative $(\varphi,\Gamma)$-module sheaf of $\scv_{0}^{\ast}$
over $\scc$, so $\scd^{\ast}$ is a locally free $\srr$-module of
rank two with $\srr$-semilinear actions of $\Gamma$ and $\varphi$.
Inside $\scd^{\ast}$ we have the sheaf of crystalline periods defined
by $\scd_{\mathrm{crys}}^{\ast}=(\scd^{\ast})^{\Gamma=1,\varphi=\alpha^{c}}$.
According to a beautiful result of Kedlaya-Pottharst-Xiao and Liu \cite{KPX,Liu},
$\scd_{\mathrm{crys}}^{\ast}$ is a line bundle on $\scc$, and for
each noble point $x$ there is a canonical identification\[
\scd_{\mathrm{crys},x}^{\ast}=\mathbf{D}_{\mathrm{crys}}(\scv_{0,x}^{\ast})^{\varphi=\alpha_{x}^{c}}.\]

\begin{thm}
\label{reg}There is a canonical morphism of sheaves of $\Lambda$-modules
\[
\mathrm{L}\mathrm{og}:\sch_{\mathrm{Iw}}^{1}(\mathbf{Q}_{p},\scv_{0})\to\sch\mathrm{om}_{\sco_{\scc}}(\scd_{\mathrm{crys}}^{\ast},\scg)\]
which interpolates the Perrin-Riou regulator map on all noble points:
given any affinoid $U\subset\scc$, any $\mathbf{z}\in\sch_{\mathrm{Iw}}^{1}(\mathbf{Q}_{p},\scv_{0})(U)$,
and any $v\in\scd_{\mathrm{crys}}^{\ast}(U)$ we have\[
\mathrm{Log}(\mathbf{z})(v)_{x}=v_{x}^{\ast}\mathrm{log}_{x}(\mathbf{z}_{x})\in\sco(\scx)\otimes_{\mathbf{Q}_{p}}E_{x}\]
for any noble point $x\in U(\overline{\mathbf{Q}_{p}})$, where \[
\mathrm{log}_{x}:H_{\mathrm{Iw}}^{1}(\mathbf{Q}_{p},\scv_{0,x})\to\sco(\scx)\otimes_{\mathbf{Q}_{p}}\mathbf{D}_{\mathrm{crys}}(\scv_{0,x})\]
is the Perrin-Riou regulator map and $v_{x}^{\ast}$ is the linear
functional induced by $v_{x}$ via the natural perfect pairing $\mathbf{D}_{\mathrm{crys}}(\scv_{0,x})\times\scd_{\mathrm{crys},x}^{\ast}\to E_{x}$.
Furthermore, $\mathrm{Log}(\mathbf{z})(v)_{x}$ has growth of order
at most $v_{p}(\alpha_{x})$ for any $x\in U(\overline{\mathbf{Q}_{p}})$.
\end{thm}
Now choose a pair of characters $\nu^{+},\nu^{-}$ as in §1.1, with
conductors $A^{\pm}$ and with $\nu^{\pm}(-1)=\pm1$, and set\begin{eqnarray*}
\mathbf{L}_{\nu^{+},\nu^{-}} & = & e^{+}\prod_{\ell|\frac{A^{-}}{(A^{+},A^{-})}}P_{\ell}(\ell^{-1}\sigma_{\ell}^{-1})\cdot\mathrm{Log}(\mathrm{res}_{p}\mathfrak{Z}_{1}(1,\nu^{+}))\\
 &  & +e^{-}\prod_{\ell|\frac{A^{+}}{(A^{+},A^{-})}}P_{\ell}(\ell^{-1}\sigma_{\ell}^{-1})\cdot\mathrm{Log}(\mathrm{res}_{p}\mathfrak{Z}_{1}(1,\nu^{-}))\end{eqnarray*}
with $P_{\ell}$ as in Theorem 1.1. (Here $e^{\pm}=\frac{1\pm\chi_{\scw}(-1)\cdot\sigma_{-1}}{2}\in\Lambda$.)
\begin{thm}
\label{twovarlfn}For any noble point $x\in U(\overline{\mathbf{Q}_{p}})\subset\scc(\overline{\mathbf{Q}_{p}})$
and any section $v$ of $\scd_{\mathrm{crys}}^{\ast}(U)$, nonvanishing
at $x$, we have\[
\mathbf{L}_{\nu^{+},\nu^{-}}(v)_{x}=C_{\nu^{+},\nu^{-}}(v_{x})\prod_{\ell|A^{+}A^{-}}P_{\ell}(\ell^{-1}\sigma_{\ell}^{-1})_{x}\cdot L_{p,\mathrm{alg}}(\mathbf{f}_{x}^{c})\in\sco(\scx)\otimes_{\mathbf{Q}_{p}}E_{x}\]
where $C_{\nu^{+},\nu^{-}}(v_{x})=e^{+}C_{\nu^{+}}(v_{x})+e^{-}C_{\nu^{-}}(v_{x})$
with $C_{\nu^{\pm}}(v_{x})$ a nonzero multiple of $L_{\mathrm{alg}}(k-1,f^{c}\otimes\nu^{-1})$.
\end{thm}
The finite product here has the effect of {}``deleting the Euler
factors'' of $L_{p,\mathrm{alg}}(\mathbf{f}_{x})$ at primes dividing
$A^{+}A^{-}$.

\subsection{A conjecture}

With an eye towards an optimal version of Theorem \ref{twovarlfn},
we offer the following conjecture.
\begin{conjecture}
There is a canonical morphism of abelian sheaves\[
\mathbf{z}_{\scc}:\scv\to\sch_{\mathrm{Iw}}^{1}(\mathbf{Q},\scv)\]
such that for any noble point $x\in\scc(\overline{\mathbf{Q}_{p}})$
and any $\gamma\in U\subset\scc$ we have\[
\mathbf{z}_{\scc}(\gamma)_{x}=\mathbf{z}_{\gamma_{x},\mathbf{f}_{x}}^{(p)}(k_{x})\in H_{\mathrm{Iw}}^{1}(\mathbf{Q},\scv_{x}),\]
where $\mathbf{z}_{\bullet,\mathbf{f}}^{(p)}$ denotes Kato's {}``optimal
zeta element'' \begin{eqnarray*}
V_{\mathbf{f}} & \to & H_{\mathrm{Iw}}^{1}(\mathbf{Q},V_{\mathbf{f}})\\
\gamma & \mapsto & \mathbf{z}_{\gamma,\mathbf{f}}^{(p)}\end{eqnarray*}
associated with a nobly refined form $\mathbf{f}$.
\end{conjecture}
It seems likely that $\mathbf{z}_{\scc}$ can be obtained by patching
suitable linear combinations of the classes $_{c,d}\mathfrak{Z}_{N,1}(\mathfrak{U},0,a(A))_{h}$
defined in §3.2 below. We shall give some partial results towards
this conjecture in \cite{Htwovar}.

\subsection*{Notation}

We write $K_{0}(N)$, $K_{1}(N)$ and $K(N)$ for the usual adelic
congruence subgroups of $\mathrm{GL}_{2}(\widehat{\mathbf{Z}})$,
and $\Gamma_{0}(N),\Gamma_{1}(N),\Gamma(N)$ for their intersections
with $\mathrm{SL}_{2}(\mathbf{Z})$. If $N=p^{n}$ is a prime power
we sometimes conflate the former groups with their projections under
the natural map \[
\mathrm{pr}_{p}:\mathrm{GL}_{2}(\widehat{\mathbf{Z}})\to\mathrm{GL}_{2}(\mathbf{Z}_{p}),\]
since $K=\mathrm{pr}_{p}^{-1}(\mathrm{pr}_{p}(K))$ in this case.
Let $I=K_{0}(p)\subset\mathrm{GL}_{2}(\mathbf{Z}_{p})$ be the subgroup
consisting of matrices which are upper-triangular modulo $p$, and
let $\Delta\subset\mathrm{GL}_{2}(\mathbf{Q}_{p})$ be the multiplicative
monoid given by\[
\Delta=\left\{ \left(\begin{array}{cc}
a & b\\
c & d\end{array}\right)\in M_{2}(\mathbf{Z}_{p})\mid\det\neq0,\, c\in p\mathbf{Z}_{p}\,\mathrm{and}\, a\in\mathbf{Z}_{p}^{\times}\right\} .\]
Note that $\Delta$ is generated by $I$ together with the matrix
$\mathrm{diag}(1,p)$.

Given $K\subset\mathrm{GL}_{2}(\widehat{\mathbf{Z}})$ open of finite
index, we define the associated open modular curve first as a complex
analytic space by \[
Y(K)(\mathbf{C})=\mathrm{GL}_{2}^{+}(\mathbf{Q})\backslash(\mathfrak{h}\times\mathrm{GL}_{2}(\mathbf{A}_{f}))/K.\]
This is a possibly disconnected Riemannian orbifold, and the determinant
map \begin{eqnarray*}
Y(K)(\mathbf{C}) & \to & \mathbf{Q}_{>0}^{\times}\backslash\mathbf{A}_{f}^{\times}/\det(K)^{\times}\simeq\widehat{\mathbf{Z}}^{\times}/\det(K)^{\times}\\
\mathrm{GL}_{2}^{+}(\mathbf{Q})(z,g_{f})K & \to & \mathbf{Q}_{>0}^{\times}\det g_{f}\det(K)^{\times}\end{eqnarray*}
induces a bijection between the set of connected components of $Y(K)(\mathbf{C})$
and the finite group $\widehat{\mathbf{Z}}^{\times}/\det(K)^{\times}$.
There is a canonical affine curve $Y(K)$ defined over $\mathbf{Q}$
such that $Y(K)(\mathbf{C})$ is the analytic space associated with
$Y(K)\times_{\mathrm{Spec}\mathbf{Q}}\mathrm{Spec}\mathbf{C}$.

We write $L,M$ for a pair of positive integers with $L+M\geq5$.
Let $Y(L,M)$ be the modular curve over $\mathbf{Q}$ representing
the functor which sends a $\mathbf{Q}$-scheme $S$ to the set of
isomorphism classes of triples $(E,e_{1},e_{2})$ where $E/S$ is
an elliptic curve and $e_{1},e_{2}\in E(S)$ are two sections of $E\to S$
such that $Le_{1}=Me_{2}=0$ and the map \begin{eqnarray*}
(\mathbf{Z}/L\mathbf{Z})\times(\mathbf{Z}/M\mathbf{Z}) & \to & E(S)\\
(a,b) & \mapsto & ae_{1}+be_{2}\end{eqnarray*}
is an injective group homomorphism. If $L|L'$ and $M|M'$ there is
a natural covering map \begin{eqnarray*}
Y(L',M') & \to & Y(L,M)\\
(E,e_{1},e_{2}) & \mapsto & (E,\tfrac{L'}{L}e_{1},\tfrac{M'}{M}e_{2}).\end{eqnarray*}
For $(a,b)\in$ $(\mathbf{Z}/L\mathbf{Z})^{\times}\times(\mathbf{Z}/M\mathbf{Z})^{\times}$,
let $\left\langle a\mid b\right\rangle $ denote the automorphism
of $Y(L,M)$ defined by $(E,e_{1},e_{2})\mapsto(E,ae_{1},abe)$, and
let $\left\langle a\mid b\right\rangle ^{\ast}$ denote the induced
pullback map on functions. Let $T_{n}$ and $T_{n}'$ denote the Hecke
operators and dual Hecke operators as defined by Kato, so $T_{\ell}=T_{\ell}'\left\langle 1/\ell\mid\ell\right\rangle $
if $\ell\nmid M$. We set $Y_{1}(M)=Y(1,M)$ as usual, and write $\left\langle 1\mid b\right\rangle =\left\langle b\right\rangle $
in this case (but not in general). Following Kato we write $Y_{1}(N)\otimes\mathbf{Q}(\zeta_{m})$
for $Y_{1}(N)\times_{\mathrm{Spec}\mathbf{Q}}\mu_{m}^{\circ}$, where
$\mu_{m}^{\circ}$ denotes the scheme of primitive $m$th roots of
unity.

We often denote {}``$\mathbf{Z}_{p}$-integral modules'' with a
superscripted $\circ$, and denote the outcome of applying $(-)\otimes_{\mathbf{Z}}\mathbf{Q}$
by removing the $\circ$. Sometimes we apply this convention in reverse:
in particular, if $M$ is a $\mathbf{Q}_{p}$-Banach module, we denote
its unit ball by $M^{\circ}$. If $X$ is a reduced rigid analytic
space, we write $\sco(X)$ for the ring of global functions on $X$,
$\sco(X)^{b}$ for the ring of bounded global functions, and $\sco(X)^{\circ}$
for the ring of power-bounded functions.

All Galois cohomology groups are taken to be continuous cohomology.

\subsection*{Acknowledgements}

The overwhelming debt this paper owes to the work of Kato will be
evident to the reader. I'd also like to acknowledge the galvanic effect
of Colmez's wonderful Bourbaki article on Kato's work \cite{ColKato};
while reading §1.5 of \cite{ColKato}, I was gripped with the conviction
that something like Theorem \ref{chern} must be true, and the present paper
grew out of that conviction. I've also benefited greatly from reading
various papers of Bellaïche, Berger, Colmez, Lei-Loeffler-Zerbes,
and Pottharst.

It's a particular pleasure to thank Barry Mazur for a number of inspiring
discussions in the early stages of this project, for his generous and invaluable encouragement, and for
providing financial support during the summer of 2013. 

Portions of this work were carried out at Boston College, l\textbf{'}Institut
de Mathématiques de Jussieu, and Columbia University; it's a pleasure to acknowledge the hospitality
of these institutions. The research leading to these results has received
funding from the European Research Council under the European Community's
Seventh Framework Programme (FP7/2007-2013) / ERC Grant agreement
n\textdegree{} 290766 (AAMOT).

\section{Background on overconvergent cohomology}

In this section we prove some foundational results on {}``étale overconvergent
cohomology.'' As noted in the introduction, we work with small formal
opens in the weight space, as opposed to the more familiar affinoid
opens. However, there is a serious payoff for this slight complication:
the filtrations which this point of view allows are so well-behaved
that they completely obliterate any possible difficulties involving
higher derived functors of inverse limits. The basic constructions
and finiteness results in this section generalize with very little
effort to yield a workable definition of étale overconvergent cohomology
for any Shimura variety: this should be particularly attractive in
the context of quaternionic Shimura curves and Kottwitz's simple Shimura
varieties, where the middle-dimensional cohomology gives the {}``right''
Galois representation.

\subsection{Locally analytic distributions}

Let $s$ be a nonnegative integer. Consider the ring of functions\[
\mathbf{A}^{s}=\left\{ f:\mathbf{Z}_{p}\to\mathbf{Q}_{p}\mid f\,\mathrm{analytic\, on\, each\,}p^{s}\mathbf{Z}_{p}-\mathrm{coset}\right\} .\]
Recall that by a fundamental result of Amice, the functions $e_{j}^{s}(x)=\left\lfloor p^{-s}j\right\rfloor !\left(\begin{array}{c}
x\\
j\end{array}\right)$ define an orthonormal basis of $\mathbf{A}^{s}$. The ring $\mathbf{A}^{s}$
is affinoid, and we set $\mathbf{B}_{s}=\mathrm{Sp}(\mathbf{A}^{s})$,
so e.g.\[
\mathbf{B}_{s}(\mathbf{C}_{p})=\left\{ x\in\mathbf{C}_{p}\mid\inf_{a\in\mathbf{Z}_{p}}|x-a|\leq p^{-s}\right\} .\]

Let $\scw$ be the rigid generic fiber of the formal scheme $\mathfrak{W}=\mathrm{Spf}(\mathbf{Z}_{p}[[\mathbf{Z}_{p}^{\times}]])$,
and let $\chi_{\scw}:\mathbf{Z}_{p}^{\times}\to\sco(\scw)^{\times}$
be the universal character induced by the inclusion of grouplike elements
$\mathbf{Z}_{p}^{\times}\subset\mathbf{Z}_{p}[[\mathbf{Z}_{p}^{\times}]]^{\times}\subset\sco(\scw)^{\times}$.
The canonical splitting \begin{eqnarray*}
\mathbf{Z}_{p}^{\times} & = & \mu_{p-1}\times(1+p\mathbf{Z}_{p})^{\times}\\
z & \mapsto & (\omega(z),\left\langle z\right\rangle )\end{eqnarray*}
induces natural identifications $\sco(\mathfrak{W})=\mathbf{Z}_{p}[\mu_{p-1}]\otimes_{\mathbf{Z}_{p}}\mathbf{Z}_{p}[[T]]$
and $\sco(\scw)=\mathbf{Z}_{p}[\mu_{p-1}]\otimes_{\mathbf{Z}_{p}}\mathbf{Q}_{p}\left\langle \left\langle T\right\rangle \right\rangle $,
where $\mathbf{Q}_{p}\left\langle \left\langle T\right\rangle \right\rangle \subset\mathbf{Q}_{p}[[T]]$
is the ring of power series convergent on the open unit disk, by mapping
the grouplike element $[1+p]\in\mathbf{Z}_{p}[[\mathbf{Z}_{p}^{\times}]]$
to $1+T$. The universal character is given by the convergent series\[
\chi_{\scw}(z)=[\omega(z)]\otimes\sum_{n=0}^{\infty}T^{n}\left(\begin{array}{c}
\frac{\log_{p}\left\langle z\right\rangle }{\log_{p}(1+p)}\\
n\end{array}\right)\in\sco(\scw),\, z\in\mathbf{Z}_{p}^{\times}.\]
The pair $(\scw,\chi_{\scw})$ is universal for pairs $(\Omega,\chi_{\Omega})$,
where $\Omega=\mathrm{Sp}A$ is an affinoid space and $\chi_{\Omega}:\mathbf{Z}_{p}^{\times}\to\sco(\Omega)^{\times}=A^{\times}$
is a continuous character; given any such pair, there is a unique
morphism $f:\Omega\to\scw$ such that $\chi_{\Omega}=f^{\ast}\chi_{\scw}$,
with $f$ characterized by $f^{\ast}[\zeta_{p-1}]=\chi_{\Omega}(\zeta_{p-1})$
and $f^{\ast}T=\chi_{\Omega}(1+p)-1$. Clearly $\chi_{\scw}$ factors
through the universal character $\chi_{\mathfrak{W}}:\mathbf{Z}_{p}^{\times}\to\sco(\mathfrak{W})^{\times}$,
and $(\mathfrak{W},\chi_{\mathfrak{W}})$ satisfies an analogous universal
property for formal schemes. If $\lambda\in\scw(\overline{\mathbf{Q}_{p}})$
is a point, we shall always denote the associated character by $\lambda$
as well. We embed $\mathbf{Z}$ into $\scw(\mathbf{Q}_{p})$ (and
into $\mathfrak{W}(\mathbf{Z}_{p})$) by mapping $k$ to the character
$\lambda_{k}(z)=z^{k-2}$.

Now let $\mathfrak{U}\subset\mathfrak{W}$ be a formal subscheme of
the form $\mathfrak{U}=\mathrm{Spf}(R_{\mathfrak{U}})$ with $R_{\mathfrak{U}}$
a $\mathbf{Z}_{p}$-flat, normal and module-finite $\mathbf{Z}_{p}[[X_{1},\dots,X_{d}]]$-algebra,
such that the induced map $\mathfrak{U}^{\mathrm{rig}}\to\scw$ on
Berthelot generic fibers is an open immersion with image contained
in an admissible affinoid open subset of $\scw$. We write $\sco(\mathfrak{U})$
interchangeably for $R_{\mathfrak{U}}$. Note that $\sco(\mathfrak{U})$
is reduced, and that we have identifications $\sco(\mathfrak{U}^{\mathrm{rig}})^{b}=\sco(\mathfrak{U})[\frac{1}{p}]$
and $\sco(\mathfrak{U}^{\mathrm{rig}})^{\circ}=\sco(\mathfrak{U})$.
We may choose an ideal $\mathfrak{a}\subset R_{\mathfrak{U}}$ containing
$p$ such that $R_{\mathfrak{U}}$ is $\mathfrak{a}$-adically separated
and complete and each $R_{\mathfrak{U}}/\mathfrak{a}^{n}$ is a finite
abelian group. Let $s_{\mathfrak{U}}$ be the least nonnegative integer
such that $\chi_{\Omega}(1+p^{s_{\mathfrak{U}}+1}x)\in\sco(\Omega)\left\langle x\right\rangle $
for some admissible affinoid $\scw\supset\Omega\Supset\mathfrak{U}^{\mathrm{rig}}$.
Define\begin{eqnarray*}
\mathcal{A}_{\mathfrak{U}}^{s,\circ} & = & \sco(\mathbf{B}_{s}\times\mathfrak{U}^{\mathrm{rig}})^{\circ}\\
 & = & \mathbf{A}^{s,\circ}\widehat{\otimes}R_{\mathfrak{U}}\end{eqnarray*}
where the indicated completion is the $1\otimes\mathfrak{a}$-adic
completion. In the Amice basis we have\[
\mathcal{A}_{\mathfrak{U}}^{s,\circ}=\left\{ \sum_{j\geq0}r_{j}e_{j}^{s}(x)\mid r_{j}\in R_{\mathfrak{U}}\,\mathrm{with}\, r_{j}\to0\,\mathfrak{a}-\mathrm{adically}\right\} .\]
The formula\[
(\gamma\cdot_{\mathfrak{U}}f)(x)=\chi_{\mathfrak{U}}(a+cx)f\left(\frac{b+dx}{a+cx}\right)\]
defines a continuous left action of $\Delta$ on $\mathcal{A}_{\mathfrak{U}}^{s,\circ}$
for any $s\geq s_{\mathfrak{U}}$. Set $\mathcal{D}_{\mathfrak{U}}^{s,\circ}=\mathrm{Hom}_{\sco(\mathfrak{U})}(\mathcal{A}_{\mathfrak{U}}^{s,\circ},\sco(\mathfrak{U}))$,
with the dual right action, and set \begin{eqnarray*}
\mathcal{A}_{\mathfrak{U}}^{s} & = & \mathcal{A}_{\mathfrak{U}}^{s,\circ}\otimes_{\mathbf{Z}}\mathbf{Q}\\
 & \cong & \sco(\mathbf{B}_{s}\times\mathfrak{U}^{\mathrm{rig}})^{b}\end{eqnarray*}
and\begin{eqnarray*}
\mathcal{D}_{\mathfrak{U}}^{s} & = & \mathcal{D}_{\mathfrak{U}}^{s,\circ}\otimes_{\mathbf{Z}}\mathbf{Q}\\
 & \cong & \mathcal{L}_{\sco(\mathfrak{U}^{\mathrm{rig}})^{b}}\left(\mathcal{A}_{\mathfrak{U}}^{s},\sco(\mathfrak{U}^{\mathrm{rig}})^{b}\right).\end{eqnarray*}
The following lemma quantifies the continuity of the $\Delta$-action
on these modules.
\begin{lem}
\label{contbasic}For any $s\geq s_{\mathfrak{U}}$ and $n\ge1$,
the principal congruence subgroup $K(p^{s+n})\subset\mathrm{GL}_{2}(\mathbf{Z}_{p})$
acts trivially on $\mathcal{A}_{\mathfrak{U}}^{s,\circ}/p^{n}$ and
$\mathcal{D}_{\mathfrak{U}}^{s,\circ}/p^{n}$.
\end{lem}
\emph{Proof. }Let $\mathcal{C}^{s}$ denote the Banach space of functions
$F:I\to\mathbf{Q}_{p}$ which are analytic on each coset of $K(p^{s})$,
in the sense that \[
F\left(\gamma\left(\begin{array}{cc}
1+p^{s}X_{1} & p^{s}X_{2}\\
p^{s}X_{3} & 1+p^{s}X_{4}\end{array}\right)\right)\in\mathbf{Q}_{p}\left\langle X_{1},X_{2},X_{3},X_{4}\right\rangle \]
for any fixed $\gamma\in I$. Regard $\mathcal{C}^{s}$ as a left
$I$-module via right translation. The key observation is that $K(p^{n+s})$
acts trivially on $\mathcal{C}^{s,\circ}/p^{n}$ and thus on $(\widehat{\mathcal{C}^{s,\circ}\otimes_{\mathbf{Z}_{p}}\sco(\mathfrak{U})})/p^{n}$.
To see this, fix $\gamma\in I$ and $F\in\mathcal{C}^{s,\circ}$,
and for $g\in K(p^{s})$ write $g=\left(\begin{array}{cc}
1+p^{s}X_{1}(g) & p^{s}X_{2}(g)\\
p^{s}X_{3}(g) & 1+p^{s}X_{4}(g)\end{array}\right)$. Define $\varphi_{F,\gamma}\in\mathbf{Z}_{p}\left\langle X_{1},\dots,X_{4}\right\rangle $
by\[
\varphi_{F,\gamma}(X_{1}(g),\dots,X_{4}(g))=F(\gamma g).\]
Then for any $h\in K(p^{s+n})$ we have\[
\mathrm{sup}_{i}|X_{i}(g)-X_{i}(gh)|\leq p^{-n},\]
so\begin{eqnarray*}
|F(\gamma g)-F(\gamma gh)| & = & |\varphi_{F,\gamma}(X_{1}(g),\dots,X_{4}(g))-\varphi_{F,\gamma}(X_{1}(gh),\dots,X_{4}(gh))|\\
 & \leq & p^{-n}\end{eqnarray*}
by Proposition 7.2.1/1 of \cite{BGR}. 

Let $B^{-}$ denote the lower-triangular matrices in $I$. A simple
calculation shows that the map $F\mapsto F\left(\begin{array}{cc}
1 & x\\
 & 1\end{array}\right)$ gives an $I$-equivariant isomorphism onto $\mathcal{A}_{\mathfrak{U}}^{s,\circ}$
from the subspace of functions $F\in\widehat{\mathcal{C}^{s,\circ}\otimes_{\mathbf{Z}_{p}}\sco(\mathfrak{U})}$
satisfying\[
F\left(\left(\begin{array}{cc}
a\\
c & d\end{array}\right)g\right)=\chi_{\mathfrak{U}}(a)F(g)\,\forall\left(\begin{array}{cc}
a\\
c & d\end{array}\right)\in B^{-}.\]
This is a closed subspace of $\widehat{\mathcal{C}^{s,\circ}\otimes_{\mathbf{Z}_{p}}\sco(\mathfrak{U})}$,
and so $\mathcal{A}_{\mathfrak{U}}^{s,\circ}/p^{n}$ injects $I$-equivariantly
into $(\widehat{\mathcal{C}^{s,\circ}\otimes_{\mathbf{Z}_{p}}\sco(\mathfrak{U})})/p^{n}$.
$\square$

Next we construct the filtration on $\mathcal{D}_{\mathfrak{U}}^{s,\circ}$
described in the introduction. For any $s>s_{\mathfrak{U}}$, define
$Q_{\mathfrak{U}}^{s,i}$ as the image of $\mathcal{D}_{\mathfrak{U}}^{s,\circ}$
in $\mathcal{D}_{\mathfrak{U}}^{s-1,\circ}/p^{i}$. Then $Q_{\mathfrak{U}}^{s,i}$
is a \emph{finitely} \emph{generated} $R_{\mathfrak{U}}/p^{i}$-submodule
of $\mathcal{D}_{\mathfrak{U}}^{s-1,\circ}/p^{i}$, stable under $\Delta$,
and the natural map\[
\mathcal{D}_{\mathfrak{U}}^{s,\circ}\to\lim_{\leftarrow i}Q_{\mathfrak{U}}^{s,i}\]
is an isomorphism. By Lemma \ref{contbasic}, $K(p^{s+i})$ acts trivially
on $Q_{\mathfrak{U}}^{s,i}$.
\begin{defn}
\textbf{\label{goodfilt}}\emph{We define\[
\mathrm{Fil}^{i}\mathcal{D}_{\mathfrak{U}}^{s,\circ}\]
as the kernel of the composite map\[
\mathcal{D}_{\mathfrak{U}}^{s,\circ}\to Q_{\mathfrak{U}}^{s,i}\to Q_{\mathfrak{U}}^{s,i}\otimes_{R_{\mathfrak{U}}/p^{i}}R_{\mathfrak{U}}/\mathfrak{a}^{i}.\]
}
\end{defn}
To unwind this, note the identification \begin{eqnarray*}
\mathcal{D}_{\mathfrak{U}}^{s,\circ} & \cong & \prod_{j\geq0}R_{\mathfrak{U}}\\
\mu & \mapsto & (\mu(e_{j}^{s}))_{j\geq0}.\end{eqnarray*}
Setting $c_{j}^{s}=\frac{\left\lfloor p^{1-s}j\right\rfloor !}{\left\lfloor p^{-s}j\right\rfloor !}$,
we have $e_{j}^{s-1}=c_{j}^{s}e_{j}^{s}$, and the coefficient $c_{j}^{s}$
is a $p$-adic integer tending uniformly to zero as $j\to\infty$,
in fact with $v_{p}(c_{j}^{s})=\left\lfloor p^{-s}j\right\rfloor $.
Therefore we obtain\[
Q_{\mathfrak{U}}^{s,i}=\bigoplus_{j\,\mathrm{with}\, v_{p}(c_{j}^{s})<i}R_{\mathfrak{U}}/p^{i-v_{p}(c_{j}^{s})}\]
with only finitely many $j$'s appearing, and thus\[
\mathcal{D}_{\mathfrak{U}}^{s,\circ}/\mathrm{Fil}^{i}\mathcal{D}_{\mathfrak{U}}^{s,\circ}=\bigoplus_{j\,\mathrm{with}\, v_{p}(c_{j}^{s})<i}R_{\mathfrak{U}}/(\mathfrak{a}^{i},p^{i-v_{p}(c_{j}^{s})}).\]
Taking the inverse limit over $i$, we easily see that $\mathcal{D}_{\mathfrak{U}}^{s,\circ}$
is separated and complete for the filtration defined by the $\mathrm{Fil}^{i}$'s,
and that this filtration possesses all the other claimed properties.

Now suppose $k\in\mathbf{Z}_{\geq2}$ and $\mathfrak{U}$ is such
that $k\in\mathfrak{U}(\mathbf{Q}_{p})$; let $\mathfrak{p}_{k}\subset\sco(\mathfrak{U})$
be the prime cutting out $k$. Let $\mathcal{A}_{k}^{s}$ be $\mathbf{A}^{s}$
with the left action\[
(\gamma\cdot_{k}f)(x)=(a+cx)^{k-2}f\left(\frac{b+dx}{a+cx}\right),\]
and let $\mathcal{D}_{k}^{s}$ be the continuous $\mathbf{Q}_{p}$-linear
dual of $\mathcal{A}_{k}^{s}$, with unit ball $\mathcal{D}_{k}^{s,\circ}$.
There is a natural $\Delta$-equivariant isomorphism\[
\mathcal{D}_{\mathfrak{U}}^{s,\circ}\otimes_{\sco(\mathfrak{U})}\sco(\mathfrak{U})/\mathfrak{p}_{k}\cong\mathcal{D}_{k}^{s,\circ}\]
which induces a $\Delta$-equivariant surjection\[
\sigma_{k}:\mathcal{D}_{\mathfrak{U}}^{s,\circ}\twoheadrightarrow\mathcal{D}_{k}^{s,\circ}.\]

Let $\scl_{k}(A)$ denote the module of polynomials $A[X]^{\mathrm{deg}\leq k-2}$
endowed with the right $\mathrm{GL}_{2}(A)$-action\[
(p\cdot_{k}\gamma)(X)=(d+cX)^{k-2}p\left(\frac{b+aX}{d+cX}\right),\]
and set in particular $\scl_{k}=\scl_{k}(\mathbf{Q}_{p})$ and $\scl_{k}^{\circ}=\scl_{k}(\mathbf{Z}_{p})$.
By a simple calculation, the map\begin{eqnarray*}
\rho_{k}:\mathcal{D}_{k}^{s,\circ} & \to & \scl_{k}^{\circ}\\
\mu & \mapsto & \int(X+x)^{k-2}\mu(x)\\
 &  & =\sum_{j=0}^{k-2}\left(\begin{array}{c}
k-2\\
j\end{array}\right)\mu(x^{j})X^{k-2-j}\end{eqnarray*}
is $\Delta$-equivariant.
\begin{defn}
\emph{The }integration map in weight $k$\emph{ is the $\Delta$-equivariant
map $i_{k}:\mathcal{D}_{\mathfrak{U}}^{s,\circ}\to\scl_{k}^{\circ}$
defined by $i_{k}=\rho_{k}\circ\sigma_{k}$.}\end{defn}
\begin{lem}
The map $i_{k}$ fits into a $\Delta$-equivariant commutative diagram\[
\xymatrix{\mathcal{D}_{\mathfrak{U}}^{s,\circ}\ar[d]\ar[r]^{i_{k}} & \scl_{k}^{\circ}\ar[d]\\
\mathcal{D}_{\mathfrak{U}}^{s,\circ}/\mathrm{Fil}^{i}\ar[r] & \scl_{k}^{\circ}/p^{i}}
\]
compatibly with varying $i$, $s$ and $\mathfrak{U}$ in the evident
manner.
\end{lem}
\emph{Proof. }We observe that $i_{k}$ can be realized as the composite
map\[
\mathcal{D}_{\mathfrak{U}}^{s,\circ}\to\mathcal{D}_{k}^{s,\circ}\to\mathcal{D}_{k}^{s-1,\circ}\overset{\rho_{k}}{\to}\scl_{k}^{\circ},\]
and the composition of the first two of these maps carries $\mathrm{Fil}^{i}\mathcal{D}_{\mathfrak{U}}^{s,\circ}$
into $p^{i}\mathcal{D}_{k}^{s-1,\circ}\subset\mathcal{D}_{k}^{s-1,\circ}$.
To see the latter, note that \begin{eqnarray*}
\mathcal{D}_{\mathfrak{U}}^{s,\circ}/\mathrm{Fil}^{i}\otimes_{R_{\mathfrak{U}}}R_{\mathfrak{U}}/\mathfrak{p}_{k} & = & Q_{\mathfrak{U}}^{s,i}\otimes_{R_{\mathfrak{U}}}R_{\mathfrak{U}}/\mathfrak{p}_{k}\\
 & \subset & (\mathcal{D}_{\mathfrak{U}}^{s-1,\circ}/p^{i})\otimes_{R_{\mathfrak{U}}}R_{\mathfrak{U}}/\mathfrak{p}_{k}\\
 & = & \mathcal{D}_{k}^{s-1,\circ}/p^{i},\end{eqnarray*}
where the first line follows from the fact that $\overline{\mathfrak{a}}\subset(p)$
in $R_{\mathfrak{U}}/\mathfrak{p}_{k}=\mathbf{Z}_{p}$. $ $ $\square$

The Dirac distribution\[
\mu_{\mathrm{Dir}}\in\mathcal{D}_{\mathfrak{U}}^{s,\circ}\]
defined by\[
\mu_{\mathrm{Dir}}(f)=f(0)\]
plays a key role in our constructions. Note that the various $\mu_{\mathrm{Dir}}$'s
are carried to each other under the {}``change of $s$'' or {}``change
of $\mathfrak{U}$'' maps, so we are justified in omitting $\mathfrak{U}$
and $s$ from the notation. We also write $\mu_{\mathrm{Dir}}$ for
the image of $\mu_{\mathrm{Dir}}$ in $\mathcal{D}_{\mathfrak{U}}^{s}$
or any quotient of $\mathcal{D}_{\mathfrak{U}}^{s,\circ}$.
\begin{lem}
\textbf{\label{compdirac}}\emph{We have $i_{k}(\mu_{\mathrm{Dir}})=X^{k-2}$.}
\end{lem}

\subsection{Étale overconvergent cohomology}

Let $Y=Y(K)$ be any modular curve with the $p$-part of $K$ contained
in $I$, so any of the modules introduced in §2.1 defines a local
system on the analytic space $Y(\mathbf{C})$ associated with $Y$.
\begin{lem}
The natural map\[
H^{1}(Y(\mathbf{C}),\mathcal{D}_{\mathfrak{U}}^{s,\circ})\to\lim_{\leftarrow i}H^{1}(Y(\mathbf{C}),\mathcal{D}_{\mathfrak{U}}^{s,\circ}/\mathrm{Fil}^{i}\mathcal{D}_{\mathfrak{U}}^{s,\circ})\]
is an isomorphism.
\end{lem}
\emph{Proof. }This map is surjective with kernel isomorphic to $\lim^{1}\left(H^{0}(Y(\mathbf{C}),\mathcal{D}_{\mathfrak{U}}^{s,\circ}/\mathrm{Fil}^{i}\mathcal{D}_{\mathfrak{U}}^{s,\circ})\right)$.
By the basic finiteness properties of $\mathcal{D}_{\mathfrak{U}}^{s,\circ}/\mathrm{Fil}^{i}\mathcal{D}_{\mathfrak{U}}^{s,\circ}$,
this is an projective system of finite abelian groups, and so its
$\lim^{1}$ vanishes. $\square$

The modules $\mathcal{D}_{\mathfrak{U}}^{s,\circ}/\mathrm{Fil}^{i}\mathcal{D}_{\mathfrak{U}}^{s,\circ}$
define locally constant sheaves of finite abelian groups on the étale
sites of $Y_{/\overline{\mathbf{Q}}}$ and $Y$. With the previous
lemma in mind, we define\[
H_{\mathrm{\acute{e}t}}^{1}(Y_{/\overline{\mathbf{Q}}},\mathcal{D}_{\mathfrak{U}}^{s,\circ})=\lim_{\leftarrow i}H_{\mathrm{\acute{e}t}}^{1}(Y_{/\overline{\mathbf{Q}}},\mathcal{D}_{\mathfrak{U}}^{s,\circ}/\mathrm{Fil}^{i}\mathcal{D}_{\mathfrak{U}}^{s,\circ}).\]
This is isomorphic to $H^{1}(Y(\mathbf{C}),\mathcal{D}_{\mathfrak{U}}^{s,\circ})$
as an $R_{\mathfrak{U}}$-module and Hecke module, and is equipped
additionally with a $p$-adically continuous%
\footnote{In the weakest possible sense: the $G_{\mathbf{Q}}$-action on $H_{\mathrm{\acute{e}t}}^{1}(Y_{/\overline{\mathbf{Q}}},\mathcal{D}_{\mathfrak{U}}^{s,\circ})\otimes_{\mathbf{Z}}\mathbf{Z}/p^{n}\mathbf{Z}$
factors through a finite quotient of $G_{\mathbf{Q}}$ for any $n$,
but please note we are \emph{not} claiming this module is separated
or complete for the $p$-adic topology!%
}$R_{\mathfrak{U}}$-linear action of $G_{\mathbf{Q}}$ commuting with
the Hecke action. If $L,M$ is a pair of integers as earlier with
$p|M$, or $N$ is an integer prime to $p$, we define in particular\begin{eqnarray*}
\mathbf{V}_{\mathfrak{U}}^{s,\circ}(L,M) & = & H_{\mathrm{\acute{e}t}}^{1}(Y(L,M)_{/\overline{\mathbf{Q}}},\mathcal{D}_{\mathfrak{U}}^{s,\circ})(2),\\
\mathbf{V}_{\mathfrak{U}}^{s,\circ}(N) & = & \mathbf{V}_{\mathfrak{U}}^{s,\circ}(1,Np).\end{eqnarray*}
We also define \[
V_{k}^{\circ}(L,M)=H_{\mathrm{\acute{e}t}}^{1}(Y(L,M)_{/\overline{\mathbf{Q}}},\scl_{k}^{\circ})(2-k).\]
Note that $\scl_{k}^{\circ}/p^{n}$ is naturally isomorphic to the
etale sheaf $\mathrm{sym}^{k-2}T_{p}E/p^{n}$. By Lemma 2.1.4, the
map $\mathcal{D}_{\mathfrak{U}}^{s,\circ}\to\scl_{k}^{\circ}/p^{j}$
factors through a map $\mathcal{D}_{\mathfrak{U}}^{s,\circ}/\mathrm{Fil}^{j}\mathcal{D}_{\mathfrak{U}}^{s,\circ}\to\scl_{k}^{\circ}/p^{j}$,
which in the inverse limit induces maps\[
\mathbf{V}_{\mathfrak{U}}^{s,\circ}(L,M)\to V_{k}^{\circ}(L,M)(k)=H_{\mathrm{\acute{e}t}}^{1}(Y(L,M)_{/\overline{\mathbf{Q}}},\scl_{k}^{\circ})(2)\]
and\[
\mathbf{V}_{\mathfrak{U}}^{s,\circ}(N)\to V_{k}^{\circ}(Np)(k).\]

Before continuing, we need some results about slope decompositions.
Recall the notion of an \emph{augmented Borel-Serre complex }of level
$K$ from \cite{HanThesis}: this is a functor $C^{\bullet}(-)=C^{\bullet}(K,-)$
from right $K$-modules to complexes of abelian groups with various
good properties, such that $H^{\ast}\left(C^{\bullet}(M)\right)\cong H^{\ast}(Y(K)(\mathbf{C}),M)$.
In particular, $C^{i}(M)=M^{\oplus r(i)}$ with $r(i)$ independent
of $M$ and zero for $i\notin[0,2]$, and $C^{\bullet}(M)$ inherits
any additional structures carried by $M$ which are compatible with
the $K$-action; furthermore, if the $K$-action on $M$ extends to
a right action of some monoid $\mathrm{S}$ with $K\subset\mathrm{S}\subset\mathrm{GL}_{2}(\mathbf{A}_{f})$,
any Hecke operator $T\in\scc_{c}^{\infty}(K\backslash\backslash\mathrm{S},\mathbf{Z})$
lifts to an endomorphism $\tilde{T}\in\mathrm{End}_{\mathrm{Ch}}(C^{\bullet}(M))$. 

Now take $K=K_{1}(Np)$, and suppose $(\mathfrak{U},h)$ is a \emph{slope
datum }in the sense of \cite{HanThesis}, i.e. we may choose an augmented
Borel-Serre complex such that $C^{\bullet}(\mathcal{D}_{\mathfrak{U}}^{s})$
admits a slope-$\leq h$ decomposition\[
C^{\bullet}(\mathcal{D}_{\mathfrak{U}}^{s})\cong C^{\bullet}(\mathcal{D}_{\mathfrak{U}}^{s})_{h}\oplus C^{\bullet}(\mathcal{D}_{\mathfrak{U}}^{s})^{h}\]
with respect to the action of $\tilde{U}_{p}$. Note that even though
$\mathcal{D}_{\mathfrak{U}}^{s}$ is not quite orthonormalizable,
all the Riesz theory arguments of \cite{BuEigen,Colemanbanachfamilies}
still go through (cf. \cite{AIS}): the determinant $\det(1-\tilde{U}_{p}X)|C^{\bullet}(\mathcal{D}_{\mathfrak{U}}^{s})$
is a well-defined Fredholm series $F_{\mathfrak{U}}(X)\in\sco(\mathfrak{U})\{\{X\}\}$
defined independently of $s$, and $C^{\bullet}(\mathcal{D}_{\mathfrak{U}}^{s})$
admits a slope-$\leq h$ decomposition if and only if $F_{\mathfrak{U}}(X)$
admits a slope-$\leq h$ factorization $F_{\mathfrak{U}}(X)=Q_{\mathfrak{U},h}(X)R(X)$,
in which case $C^{\bullet}(\mathcal{D}_{\mathfrak{U}}^{s})_{h}=\ker Q_{\mathfrak{U},h}^{\ast}(\tilde{U}_{p})$.%
\footnote{Here $Q_{\mathfrak{U},h}(X)\in\sco(\mathfrak{U})[X]$ has leading
coefficient a unit in $\sco(\mathfrak{U})[\frac{1}{p}]$, and $Q^{\ast}(X)=X^{\deg Q}Q(1/X)$.%
} $ $ Taking the cohomology of this decomposition yields a Hecke-
and Galois-stable direct sum decomposition $\mathbf{V}_{\mathfrak{U}}^{s}(N)=\mathbf{V}_{\mathfrak{U},h}^{s}(N)\oplus\mathbf{V}_{\mathfrak{U}}^{s}(N)^{h}$,
and $\mathbf{V}_{\mathfrak{U},h}^{s}(N)$ is a finite projective $\sco(\mathfrak{U})[\tfrac{1}{p}]$-module.
We define $\mathbf{V}_{\mathfrak{U},h}^{s,\circ}(N)$ as the image
of $\mathbf{V}_{\mathfrak{U}}^{s,\circ}(N)$ in $\mathbf{V}_{\mathfrak{U},h}^{s}(N)$,
so $\mathbf{V}_{\mathfrak{U},h}^{s}(N)=\mathbf{V}_{\mathfrak{U},h}^{s,\circ}(N)\otimes_{\mathbf{Z}_{p}}\mathbf{Q}_{p}$.
The following proposition is a key finiteness result.
\begin{prop}
\label{fin}The $\sco(\mathfrak{U})$-module \textbf{\emph{$\mathbf{V}_{\mathfrak{U},h}^{s,\circ}(N)$}}
is finitely presented, and in particular is separated and complete
for the $\mathfrak{a}$-adic topology, and thus for the p-adic topology
as well.
\end{prop}
\emph{Proof. }Fix the Borel-Serre complex $C^{\bullet}(-)=C^{\bullet}(Y_{1}(Np),-)$
for which $(\mathfrak{U},h)$ is a slope datum, and let $Z^{1}(-)$
denote the cocycles in degree $1$. The direct sum decomposition $C^{\bullet}(\mathcal{D}_{\mathfrak{U}}^{s})=C^{\bullet}(\mathcal{D}_{\mathfrak{U}}^{s})_{h}\oplus C^{\bullet}(\mathcal{D}_{\mathfrak{U}}^{s})^{h}$
induces an analogous decomposition of the cocycles. We have the commutative
diagram\[
\xymatrix{Z^{1}(\mathcal{D}_{\mathfrak{U}}^{s})_{h}\ar@{->>}[r]\ar[d]_{i_{h}} & H^{1}(\mathcal{D}_{\mathfrak{U}}^{s})_{h}=\mathbf{V}_{\mathfrak{U},h}^{s}(N)\ar[d]_{i_{h}}\\
Z^{1}(\mathcal{D}_{\mathfrak{U}}^{s})\ar@{->>}[r]\ar@/{}_{1pc}/[u]_{\mathrm{pr}_{h}} & H^{1}(\mathcal{D}_{\mathfrak{U}}^{s})=\mathbf{V}_{\mathfrak{U}}^{s}(N)\ar@/{}_{1pc}/[u]_{\mathrm{pr}_{h}}\\
Z^{1}(\mathcal{D}_{\mathfrak{U}}^{s,\circ})\ar@{->>}[r]\ar[u] & H^{1}(\mathcal{D}_{\mathfrak{U}}^{s,\circ})=\mathbf{V}_{\mathfrak{U}}^{s,\circ}(N)\ar[u]}
\]
where $i_{h}$ (resp. $\mathrm{pr}_{h}$) denotes the natural inclusion
(resp. projection). Choose elements $m_{1},\dots,m_{g}\in\mathbf{V}_{\mathfrak{U},h}^{s,\circ}(N)$
which generate $\mathbf{V}_{\mathfrak{U},h}^{s}(N)$ as an $\sco(\mathfrak{U})[\frac{1}{p}]$-module.
Let $M$ be the $\sco(\mathfrak{U})$-submodule of $\mathbf{V}_{\mathfrak{U},h}^{s}(N)$
generated by $m_{1},\dots,m_{g}$, and set $Q=\mathbf{V}_{\mathfrak{U},h}^{s,\circ}(N)/M$.
Clearly any element of $Q$ is killed by a finite power of $p$, and
it's enough to show that $Q$ has finite exponent $p^{e}<\infty$,
since then $\mathbf{V}_{\mathfrak{U},h}^{s,\circ}(N)\subseteq p^{-e}M$
as $\sco(\mathfrak{U})$-modules and $p^{-e}M$ is a finite $\sco(\mathfrak{U})$-module
by construction.

Suppose $Q$ does not have finite exponent, so we may choose cohomology
classes $c_{1},c_{2},\dots,c_{n},\dots\phantom{}\in\mathbf{V}_{\mathfrak{U}}^{s,\circ}(N)$
such that the image of $\mathrm{pr}_{h}(c_{n})$ in $Q$ has exponent
$p^{n}$. The module $\mathbf{V}_{\mathfrak{U},h}^{s}(N)$ is naturally
a \emph{p-}adic Banach space, and the sequence $\mathrm{pr}_{h}(c_{n})$
is unbounded. Now choose any cocycles $z_{n}\in Z^{1}(\mathcal{D}_{\mathfrak{U}}^{s,\circ})$
lifting the $c_{n}$'s. The sequence of $z_{n}$'s is clearly bounded
in the Banach topology on $Z^{1}(\mathcal{D}_{\mathfrak{U}}^{s})$,
and $\mathrm{pr}_{h}:Z^{1}(\mathcal{D}_{\mathfrak{U}}^{s})\to Z^{1}(\mathcal{D}_{\mathfrak{U}}^{s})_{h}$
is a continuous map of $p$-adic Banach spaces, so the sequence of
cocycles $\mathrm{pr}_{h}(z_{n})$ is bounded as well. But the topmost
horizontal arrow in the diagram is a continuous map of $p$-adic Banach
spaces carrying $\mathrm{pr}_{h}(z_{n})$ to $\mathrm{pr}_{h}(c_{n})$,
so we have a contradiction. $\square$

The map $\mathbf{V}_{\mathfrak{U}}^{s+1}(N)\to\mathbf{V}_{\mathfrak{U}}^{s}(N)$
induces canonical \emph{isomorphisms $\mathbf{V}_{\mathfrak{U},h}^{s+1}(N)\cong\mathbf{V}_{\mathfrak{U},h}^{s}(N)$},
and we write $\scv_{\mathfrak{U},h}(N)$ for the inverse limits along
these isomorphisms. 
\begin{prop}
The modules $\mathbf{V}_{\mathfrak{U},h}^{s,\circ}(N)$ and $\mathbf{V}_{\mathfrak{U},h}^{s}(N)$
are p-torsion-free, and their natural $G_{\mathbf{Q}}$-actions are
$ $$p$-adically continuous.
\end{prop}
\emph{Proof. }Immediate. $\square$

\subsection{Reconstructing the eigencurve}

Let $\mathbf{B}[r]=\mathrm{Sp}\mathbf{Q}_{p}\left\langle p^{r}X\right\rangle $
be the rigid disk of radius $p^{r}$, with $\mathbf{A}^{1}=\cup_{r}\mathbf{B}[r]$.
Now let $F\in\sco(\mathfrak{W})\{\{X\}\}$ be a Fredholm series, with
$\scz\subset\scw\times\mathbf{A}^{1}$ the associated Fredholm hypersurface.
For any admissible open $U\subset\scw$ and any $h\in\mathbf{Q}$,
$U\times\mathbf{B}[h]$ is admissible open in $\scw\times\mathbf{A}^{1}$,
and we define an admissible open subset $\scz_{U,h}\subset\scz$ by
$\scz_{U,h}=\scz\cap(U\times\mathbf{B}[h])$. If $U$ is affinoid
then $\scz_{U,h}$ is affinoid as well (and in fact $\scz_{U,h}=\left(\mathrm{Sp}\sco(U)\left\langle p^{h}X\right\rangle /(F(X))\right)$).
We say $\scz_{U,h}$ is \emph{slope-adapted }if the natural map $\scz_{U,h}\to U$
is finite and flat.
\begin{lem}[Buzzard]
 There is an admissible covering of $\scz$ by slope-adapted affinoids
$\scz_{U,h}$.
\end{lem}
\emph{Proof. }This is Theorem 4.6 of \cite{BuEigen}. $\square$
\begin{lem}[Buzzard]
 If $U\subset\scw$ is an affinoid subdomain and $h\in\mathbf{Q}$
with $\scz_{U,h}$ slope-adapted, we can choose an affinoid subdomain
$V\subset\scw$ and an $h'\geq h$ such that $U\Subset V$, $\scz_{V,h'}$
is slope-adapted, and $\scz_{V,h'}\cap(U\times\mathbf{A}^{1})=\scz_{U,h}$.
\end{lem}
\emph{Proof. }Immediate from Lemma 4.5 of \cite{BuEigen}. Note that
$\scz_{U,h}=\scz_{U,h'}$. $\square$
\begin{lem}
For some index set $I$, we may choose pairs of affinoid subdomains
$U_{i}'\Subset U_{i}$ and rationals $h_{i}$ such that $\scz_{U_{i},h_{i}}\cap(U_{i}'\times\mathbf{A}^{1})=\scz_{U_{i}',h_{i}}$
and both $(\scz_{U_{i},h_{i}})_{i\in I}$ and $(\scz_{U_{i}',h_{i}})_{i\in I}$
are admissible coverings of $\scz$ by slope-adapted affinoids.
\end{lem}
\emph{Proof. }Immediate from the previous two lemmas. $\square$
\begin{lem}
\label{goodcov}Notation as in the previous lemma, we may choose formal
opens $\mathfrak{U}_{i}\subset\mathfrak{W}$ such that $U_{i}'\subset\mathfrak{U}_{i}^{\mathrm{rig}}\subset U_{i}$
and such that $(\scz_{\mathfrak{U}_{i}^{\mathrm{rig}},h_{i}})_{i\in I}$
is an admissible covering of $\scz$ by slope-adapted admissible opens.
\end{lem}
\emph{Proof. }For each $i$, choose a finite map $f_{i}^{\ast}:\mathbf{Q}_{p}\left\langle X_{1},\dots,X_{d}\right\rangle \to\sco(U{}_{i})$
for some $d$ such that the associated map $f_{i}:U_{i}\to\mathrm{Sp}\mathbf{Q}_{p}\left\langle X_{1},\dots,X_{d}\right\rangle $
carries $U_{i}'$ into $\mathrm{Sp}\mathbf{Q}_{p}\left\langle p^{-r}X_{1},\dots,p^{-r}X_{d}\right\rangle $
for some $r\in\mathbf{Q}_{>0}$. The existence of such a map is the
definition of the relative compacity $U_{i}'\Subset U_{i}$. By Corollary
6.4.1/6 of \cite{BGR}, $\sco(U_{i})^{\circ}$ is module-finite over
$\mathbf{Z}_{p}\left\langle X_{1},\dots,X_{d}\right\rangle $. Now
let \[
R_{i}'=\left(\mathbf{Z}_{p}[[X_{1},\dots,X_{d}]]\otimes_{\mathbf{Z}_{p}\left\langle X_{1},\dots,X_{d}\right\rangle ,f_{i}^{\ast}}\sco(U_{i})^{\circ}\right)/p-\mathrm{power\, torsion}\]
and let $R_{i}$ be the normalization of $R_{i}'$; setting $\mathfrak{U}_{i}=\mathrm{Spf}(R_{i})$
gives a formal scheme of the desired type with $\mathfrak{U}_{i}^{\mathrm{rig}}\subset U_{i}$,
and the map $\sco(U_{i})\to\sco(U_{i}')$ clearly factors over a map
$R_{i}[\frac{1}{p}]\to\sco(U_{i}')$ (since $R_{i}'[\frac{1}{p}]\cong R_{i}[\frac{1}{p}]$
and the $X_{i}$'s are topologically nilpotent in $\sco(U_{i}')$
by design), so $U_{i}'\subset\mathfrak{U}_{i}^{\mathrm{rig}}$ as
desired. $\square$

Notation as in the previous subsection, fix a choice of augmented
Borel-Serre complex, and let $F\in\sco(\mathfrak{W})\{\{X\}\}$ be
the Fredholm series such that $F|_{\mathfrak{U}}=\det(1-\tilde{U}_{p}X)|C^{\bullet}(\mathcal{D}_{\mathfrak{U}}^{s})$
for all $\mathfrak{U}\subset\mathfrak{W}$. Let $\scz\subset\scw\times\mathbf{A}^{1}$
be the associated Fredholm hypersurface. Given any slope-adapted admissible
open $\scz_{\mathfrak{U}^{\mathrm{rig}},h}\subset\scz$, we have a
natural identification $\sco(\scz_{\mathfrak{U}^{\mathrm{rig}},h})=\sco(\mathfrak{U}^{\mathrm{rig}})[X]/Q_{\mathfrak{U},h}(X)$,
and the $\sco(\mathfrak{U}^{\mathrm{rig}})$-module \[
\scv_{\mathfrak{U}^{\mathrm{rig}},h}(N):=\scv_{\mathfrak{U},h}(N)\otimes_{\sco(\mathfrak{U})[\frac{1}{p}]}\sco(\mathfrak{U}^{\mathrm{rig}})\]
has a natural structure as an $\sco(\scz_{\mathfrak{U}^{\mathrm{rig}},h})$-module,
by letting $X$ act as $U_{p}^{-1}$. These actions are compatible
as $\mathfrak{U}$ and $h$ vary, and in particular we may glue the
$\scv_{\mathfrak{U}^{\mathrm{rig}},h}(N)$'s over a chosen admissible
covering of the type provided by Lemma \ref{goodcov} into a coherent
sheaf $\scv(N)$ over $\scz$ such that $ $$\scv(N)(\scz_{\mathfrak{U}^{\mathrm{rig}},h})\cong\scv_{\mathfrak{U}^{\mathrm{rig}},h}(N)$.
For any constituent $\scz_{\mathfrak{U}^{\mathrm{rig}},h}$ of our
chosen covering, let $\mathbf{T}_{\mathfrak{U},h}(N)$ denote the
image of $\sco(\mathfrak{U})[\frac{1}{p}]\otimes_{\mathbf{Z}}\mathbf{T}$
in $\mathrm{End}_{\sco(\mathfrak{U})[\frac{1}{p}]}\left(\scv_{\mathfrak{U},h}(N)\right)$.
This is clearly a module-finite $\sco(\mathfrak{U})[\frac{1}{p}]$-algebra,
and we have natural quasi-Stein rigid spaces $\scc_{\mathfrak{U}^{\mathrm{rig}},h}(N)$,
finite over $\scz_{\mathfrak{U}^{\mathrm{rig}},h}$, such that \[
\sco(\scc_{\mathfrak{U}^{\mathrm{rig}},h}(N))\cong\mathbf{T}_{\mathfrak{U},h}(N)\otimes_{\sco(\mathfrak{U})[\frac{1}{p}]}\sco(\mathfrak{U}^{\mathrm{rig}}).\]
We finally obtain the eigencurve $\scc(N)$ by gluing these rigid
spaces over our chosen covering of $\scz$. Since the $\scv_{\mathfrak{U}^{\mathrm{rig}},h}(N)$'s
are naturally finite $\sco(\scc_{\mathfrak{U}^{\mathrm{rig}},h}(N))$-modules,
we may glue them into a coherent sheaf $\scv(N)$ over $\scc(N)$
(under the spectral projection $\pi:\scc(N)\to\scz$ we have $\pi_{\ast}\scv(N)=\scv(N)$,
with a slight abuse of notation). The Galois action on each $\scv_{\mathfrak{U}^{\mathrm{rig}},h}(N)$
glues into an action on $\scv(N)$. Finally, for any number field
$K$ and any finite set of places $S$ of $K$ containing all places
dividing $Np\infty$, the Galois cohomology modules $H^{i}(G_{K,S},\scv_{\mathfrak{U}^{\mathrm{rig}},h}(N))$
glue (cf. \cite{PottharstAFFSSG}) into a coherent sheaf $\sch^{i}(G_{K,S},\scv(N))$
over $\scc(N)$.
\begin{prop}
The coherent sheaf $\scv(N)$ over $\scc(N)$ is torsion-free, and
$\scv(N)|_{\scc_{0}^{M-\mathrm{new}}(N)}$ has generic rank $2\tau(N/M)$,
where $\tau(A)$ denotes the number of divisors of $A$. In particular,
the pullback of $\scv(N)$ to the normalization $\scc$ of $\scc_{0}^{N-\mathrm{new}}(N)$
is locally free of rank two.
\end{prop}
\emph{Proof (sketch). }Torsion-freeness follows from the torsion-freeness
of $H^{1}(Y_{1}(Np),\mathcal{D}_{\mathfrak{U}}^{s})$ over $\sco(\mathfrak{U})[\frac{1}{p}]$,
which in turn follows from the vanishing of $H^{0}(Y_{1}(Np)(\mathbf{C}),\mathcal{D}_{\lambda}^{s})$
(here $\lambda\in\scw$ is an arbitrary weight) by examining the long
exact sequence in cohomology for\[
0\to\mathcal{D}_{\mathfrak{U}}^{s}\overset{\cdot a}{\to}\mathcal{D}_{\mathfrak{U}}^{s}\to\mathcal{D}_{\lambda}^{s}\to0,\]
where e.g. $a\in\sco(\mathfrak{U})[\frac{1}{p}]$ is a prime element
with zero locus $\lambda$. The statements about ranks follow from
an straightforward analysis of the specializations of $\scv(N)$ at
crystalline points of conductor $M$ and non-critical slope, by combining
Stevens's control theorem, the Eichler-Shimura isomorphism, and basic
newform theory. $\square$

Let $\mathbf{T}_{\mathfrak{U},h}^{s,\circ}(N)$ denote the image of
$\mathbf{T}\otimes_{\mathbf{Z}}\sco(\mathfrak{U})$ in $\mathrm{End}_{\sco(\mathfrak{U})}\left(\mathbf{V}_{\mathfrak{U},h}^{s,\circ}(N)\right)$.
The ring $\mathbf{T}_{\mathfrak{U},h}^{s,\circ}(N)$ is a $\mathbf{Z}_{p}$-flat
complete semilocal Noetherian ring of dimension two. Set $\mathfrak{C}_{\mathfrak{U},h}^{s}(N)=\mathrm{Spf}\mathbf{T}_{\mathfrak{U},h}^{s,\circ}(N)$.
Since $\mathbf{T}_{\mathfrak{U},h}^{s,\circ}(N)[\frac{1}{p}]\cong\mathbf{T}_{\mathfrak{U},h}(N)$,
the generic fiber of $\mathfrak{C}_{\mathfrak{U},h}^{s}(N)$ recovers
$\scc_{\mathfrak{U}^{\mathrm{rig}},h}(N)$. Is it possible to glue
the formal schemes $\mathfrak{C}_{\mathfrak{U},h}^{s}(N)$, and the
Galois representations $\mathbf{V}_{\mathfrak{U},h}^{s,\circ}(N)$
over them, in any meaningful way?

\section{Big zeta elements}

Throughout this section, we let $L,M$ denote a pair of positive integers
with $p|M$, and let $N$ denote a positive integer with $p\nmid N$.
The main reference for this section is §8 of \cite{Kato}. We largely
omit proofs in this section, but this seems a reasonable decision
to us: these proofs are (unsurprisingly) given by cross-pollinating
the proofs in  \cite{Kato} with Theorem \ref{chern} below.

\subsection{Chern class maps and zeta elements}
\begin{thm}
\label{chern}There is a canonical family of maps \[
\mathrm{Ch}_{L,M}^{s}(\mathfrak{U},r):\lim_{\substack{\leftarrow\\
j}
}K_{2}(Y(Lp^{j},Mp^{j}))\to H^{1}\left(G_{\mathbf{Q}},\mathbf{V}_{\mathfrak{U}}^{s,\circ}(L,M)(-r)\right)\]
with the following properties.

\textbf{\emph{i.}} The diagram\emph{\[
\xymatrix{\lim K_{2}(Y(Lp^{j},Mp^{j}))\ar[rr]^{\mathrm{Ch}_{L,M}^{s+1}(\mathfrak{U},r)\quad\quad}\ar[drr]_{\mathrm{Ch}_{L,M}^{s}(\mathfrak{U}',r)} &  & H^{1}\left(G_{\mathbf{Q}},\mathbf{V}_{\mathfrak{U}}^{s+1,\circ}(L,M)(-r)\right)\ar[d]\\
 &  & H^{1}\left(G_{\mathbf{Q}},\mathbf{V}_{\mathfrak{U}'}^{s,\circ}(L,M)(-r)\right)}
\]
}commutes for any $\mathfrak{U}'\subseteq\mathfrak{U}$ and $s\geq s_{\mathfrak{U}}$,
where the vertical arrow is induced by the map $\mathcal{D}_{\mathfrak{U}}^{s+1,\circ}\to\mathcal{D}_{\mathfrak{U}}^{s,\circ}\to\mathcal{D}_{\mathfrak{U}'}^{s,\circ}$.

\textbf{\emph{ii.}}\emph{ }For any $k\in\mathbf{Z}_{\geq2}\cap\mathfrak{U}(\mathbf{Q}_{p})$,
the diagram \emph{\[
\xymatrix{\lim K_{2}(Y(Mp^{j},Np^{j}))\ar[rr]^{\mathrm{Ch}_{L,M}^{s}(\mathfrak{U},r)\quad\quad}\ar[drr]_{\mathrm{Ch}_{L,M}(k,r,k-1)} &  & H^{1}\left(G_{\mathbf{Q}},\mathbf{V}_{\mathfrak{U}}^{s,\circ}(L,M)(-r)\right)\ar[d]\\
 &  & H^{1}\left(G_{\mathbf{Q}},V_{k}^{\circ}(L,M)(k-r)\right)}
\]
}commutes, where $\mathrm{Ch}_{L,M}(k,r,r')$ is Kato's Chern class
map and the vertical arrow is induced by the integration map in weight
$k$.
\end{thm}
\emph{Proof. }For brevity, set $Y_{j}=Y(Lp^{j},Mp^{j})$. For some
large fixed $m$, the diagram\[
\xymatrix{\underset{\leftarrow}{\lim}K_{2}(Y_{j+m})\ar[d]\\
\underset{\leftarrow}{\lim}H_{\mathrm{\acute{e}t}}^{2}\left(Y_{j+m},(\mathbf{Z}/p^{j})(2)\right)\ar[d]^{\cup X^{k-2}\otimes\zeta_{p^{j}}^{\otimes-r}}\ar[dr]^{\cup\mu_{\mathrm{Dir}}\otimes\zeta_{p^{j+m}}^{\otimes-r}}\\
\underset{\leftarrow}{\lim}H_{\mathrm{\acute{e}t}}^{2}\left(Y_{j+m},(\scl_{k}^{\circ}/p^{j})(2-r)\right)\ar[d]_{\mathrm{tr}} & \underset{\leftarrow}{\lim}H_{\mathrm{\acute{e}t}}^{2}\left(Y_{j+m},(\mathcal{D}_{\mathfrak{U}}^{s,\circ}/\mathrm{Fil}^{j+m-s}\mathcal{D}_{\mathfrak{U}}^{s,\circ})(2-r)\right)\ar[d]_{\mathrm{tr}}\\
\underset{\leftarrow}{\lim}H_{\mathrm{\acute{e}t}}^{2}\left(Y(L,M),(\scl_{k}^{\circ}/p^{j})(2-r)\right)\ar[d]^{d} & \underset{\leftarrow}{\lim}H_{\mathrm{\acute{e}t}}^{2}\left(Y(L,M),(\mathcal{D}_{\mathfrak{U}}^{s,\circ}/\mathrm{Fil}^{j+m-s}\mathcal{D}_{\mathfrak{U}}^{s,\circ})(2-r)\right)\ar[d]^{d}\\
\underset{\leftarrow}{\lim}H^{1}\left(G_{\mathbf{Q}},H_{\mathrm{\acute{e}t}}^{1}\left(Y(L,M)_{/\overline{\mathbf{Q}}},(\scl_{k}^{\circ}/p^{j})\right)(2-r)\right)\ar[d]^{\wr} & \underset{\leftarrow}{\lim}H^{1}\left(G_{\mathbf{Q}},H_{\mathrm{\acute{e}t}}^{1}\left(Y(L,M)_{/\overline{\mathbf{Q}}},\mathcal{D}_{\mathfrak{U}}^{s,\circ}/\mathrm{Fil}^{j+m-s}\mathcal{D}_{\mathfrak{U}}^{s,\circ}\right)(2-r)\right)\ar[d]^{\wr}\\
H^{1}\left(G_{\mathbf{Q}},V_{k}^{\circ}(L,M)(k-r)\right) & H^{1}\left(G_{\mathbf{Q}},\mathbf{V}_{\mathfrak{U}}^{s,\circ}(L,M)(-r)\right)\ar[l]_{i_{k}}}
\]
commutes for any $k\in\mathfrak{U}\cap\mathbf{Z}_{\geq2}$. Here all
inverse limits are taken over $j$, the topmost vertical arrow is
the étale Chern class map explained in Kato, the arrows labeled {}``tr''
are induced by the trace map on cohomology, and the arrows labeled
{}``$d$'' are the edge maps in the Hochschild-Serre spectral sequence.
Crucially, the diagonal arrow is well-defined, since by Lemma \ref{contbasic}
$K(p^{j+m})$ acts trivially on $\mathcal{D}_{\mathfrak{U}}^{s,\circ}/\mathrm{Fil}^{j+m-s}\mathcal{D}_{\mathfrak{U}}^{s,\circ}$.
Kato defines $\mathrm{Ch}_{L,M}(k,r,k-1)$ as the composite of all
the maps in the lefthand column. We simply define $\mathrm{Ch}_{L,M}^{s}(\mathfrak{U},r)$
as the composite map from the upper left to the lower right. $\square$
\begin{prop}
The map $\mathrm{Ch}_{L,M}^{s}(\mathfrak{U},r)$ satisfies the equivariance
formulas\[
\left\langle a\mid b\right\rangle \mathrm{Ch}_{L,M}^{s}(\mathfrak{U},r)=\chi_{\mathfrak{U}}(a)(ab)^{-r}\mathrm{Ch}_{L,M}^{s}(\mathfrak{U},r)\left\langle a\mid b\right\rangle \]
and\[
T_{m}'\mathrm{Ch}_{L,M}^{s}(\mathfrak{U},r)=\chi_{\mathfrak{U}}(m)\mathrm{Ch}_{L,M}^{s}(\mathfrak{U},r)T_{m}',\]
for any integers $a,b,m$ with $(m,Lp)=(a,Lp)=(b,Mp)=1$.
\end{prop}
These maps, much like adult mayflies, exist for exactly one purpose.
Let $ $\[
(_{c,d}z_{Lp^{j},Mp^{j}})_{j\geq1}\in\lim_{\substack{\leftarrow\\
j}
}K_{2}(Y_{j})\]
be Kato's norm-compatible system of zeta elements (§2 of \cite{Kato}).
\begin{defn}
\emph{We define\begin{eqnarray*}
_{c,d}\mathfrak{z}_{L,M}^{s}(\mathfrak{U},r) & = & \mathrm{Ch}_{L,M}^{s}(\mathfrak{U},r)\left((_{c,d}z_{Lp^{j},Mp^{j}})_{j\geq1}\right)\\
 & \in & H^{1}\left(G_{\mathbf{Q}},\mathbf{V}_{\mathfrak{U}}^{s,\circ}(L,M)(-r)\right)\end{eqnarray*}
for integers $c,d$ with $(c,6pL)=(d,6pM)=1$.}
\end{defn}
By Theorem \ref{chern}, we immediately deduce
\begin{prop}
\emph{For any $k\in\mathbf{Z}_{\geq2}\cap\mathfrak{U}$, we have an
equality \[
i_{k}\left(_{c,d}\mathfrak{z}_{L,M}^{s}(\mathfrak{U},r)\right)=\phantom{}_{c,d}\mathbf{z}_{L,M}^{(p)}(k,r,k-1)\]
in $H^{1}\left(G_{\mathbf{Q}},V_{k}^{\circ}(M,N)(k-r)\right)$ where
}\textup{$\phantom{}_{c,d}\mathbf{z}_{L,M}^{(p)}(k,r,k-1)$ }\emph{is
the p-adic zeta element defined in §8.4 of} \cite{Kato}\emph{.}
\end{prop}
By the argument in §8.5-8.6 of \cite{Kato}, the class $_{c,d}\mathfrak{z}_{L,M}^{s}(\mathfrak{U},r)$
is unramified away from $p$.

\subsection{Cyclotomic zeta elements}

In this subsection we \emph{p-}adically interpolate the construction
given in §8.9 of \cite{Kato}. Let $m,A$ be positive integers with
$p\nmid A$, let $a(A)$ denote a residue class modulo $A$, and let
$c,d$ be any integers with $(cd,6pmA)=(d,N)=1$. Choose any integers
$L,M$ with $mA|L$, $L|M$, $N|M$, $\mathrm{prime}(L)=\mathrm{prime}(mpA)$,
and $\mathrm{prime}(M)=\mathrm{prime}(mpAN)$, and let\[
t_{m,a(A)}:\mathbf{V}_{\mathfrak{U}}^{s,\circ}(L,M)\to\mathbf{V}_{\mathfrak{U}}^{s,\circ}(N)\otimes_{\mathbf{Z}}\mathbf{Z}[G_{\mathbf{Q}(\zeta_{m})/\mathbf{Q}}]\]
be the trace map induced by the morphism $Y(L,M)\to Y_{1}(Np)\otimes\mathbf{Q}(\zeta_{m})$
as in Kato. Define $_{c,d}\mathfrak{z}_{N,m}^{s}(\mathfrak{U},r,a(A))$
as the image of $_{c,d}\mathfrak{z}_{L,M}^{s}(\mathfrak{U},r)$ under
the homomorphism\begin{eqnarray*}
H^{1}\left(G_{\mathbf{Q}},\mathbf{V}_{\mathfrak{U}}^{s,\circ}(L,M)(-r)\right) & \overset{t_{m,a(A)}}{\to} & H^{1}\left(G_{\mathbf{Q}},\mathbf{V}_{\mathfrak{U}}^{s,\circ}(N)(-r)\otimes_{\mathbf{Z}}\mathbf{Z}[G_{\mathbf{Q}(\zeta_{m})/\mathbf{Q}}]\right)\\
 & \overset{\sim}{\to} & H^{1}\left(G_{\mathbf{Q}(\zeta_{m})},\mathbf{V}_{\mathfrak{U}}^{s,\circ}(N)(-r)\right).\end{eqnarray*}

\begin{prop}
\label{compaA}For any $k\in\mathbf{Z}_{\geq2}\cap\mathfrak{U}$,
we have an equality \[
i_{k}\left(_{c,d}\mathfrak{z}_{N,m}^{s}(\mathfrak{U},r,a(A))\right)=\phantom{}_{c,d}\mathbf{z}_{1,Np,m}^{(p)}(k,r,k-1,a(A),\mathrm{prime}(mpA))\]
in $H^{1}\left(G_{\mathbf{Q}(\zeta_{m})},V_{k}^{\circ}(Np)(k-r)\right)$,
where $\phantom{}_{c,d}\mathbf{z}_{1,Np,m}^{(p)}(k,r,k-1,a(A),\mathrm{prime}(mpA))$
is the p-adic zeta element in equation (8.1.2) of \cite{Kato}.
\end{prop}
Note that Kato defines $\phantom{}_{c,d}\mathbf{z}_{1,N,m}^{(p)}(k,r,k-1,a(A),\mathrm{prime}(mpA))$
without assuming $p\nmid A$, but we require this since $t_{m,a(A)}$
does not respect $\Gamma_{0}(p)$-structures when $p|A$.
\begin{prop}
\label{eulerbasic}Notation as above, let $\ell$ be any prime with
$(\ell,cd)=1$. Then\emph{\[
\mathrm{Cor}_{\mathbf{Q}(\zeta_{m\ell})/\mathbf{Q}(\zeta_{m})}\left(_{c,d}\mathfrak{z}_{N,\ell m}^{s}(\mathfrak{U},r,a(A))\right)=\begin{cases}
_{c,d}\mathfrak{z}_{N,m}^{s}(\mathfrak{U},r,a(A)) & \mathrm{if}\,\ell|mpA\\
P_{\ell}(\sigma_{\ell}^{-1})\cdot{}_{c,d}\mathfrak{z}_{N,m}^{s}(\mathfrak{U},r,a(A)) & \mathrm{if}\ell\nmid mpA,\end{cases}\]
}where\[
P_{\ell}(\sigma_{\ell}^{-1})=(1-T_{\ell}'\ell^{-r}\sigma_{\ell}^{-1}+\delta_{\ell\nmid N}\left\langle \ell\right\rangle ^{-1}\chi_{\scw}(\ell)\ell^{-1-2r}\sigma_{\ell}^{-2})\in\sco(\mathfrak{U})\otimes_{\mathbf{Z}}\mathbf{T}[G_{\mathbf{Q}(\zeta_{m})/\mathbf{Q}}].\]

\end{prop}
\emph{Proof. }Follows formally from the properties of all the maps
defined so far and Kato's analysis of the Euler system relations in
$K_{2}$. $\square$

Suppose now that $(\mathfrak{U},h)$ is a slope datum. Define $_{c,d}\mathfrak{z}_{N,m}^{s}(\mathfrak{U},r,a(A))_{h}$
as the image of $_{c,d}\mathfrak{z}_{N,m}^{s}(\mathfrak{U},r,a(A))$
under the natural map\[
H^{1}\left(G_{\mathbf{Q}(\zeta_{m})},\mathbf{V}_{\mathfrak{U}}^{s,\circ}(N)(-r)\right)\to H^{1}\left(G_{\mathbf{Q}(\zeta_{m})},\mathbf{V}_{\mathfrak{U},h}^{s,\circ}(N)(-r)\right).\]
Now define\begin{eqnarray*}
H_{\mathrm{Iw}}^{1}(G_{\mathbf{Q}(\zeta_{m})},\mathbf{V}_{\mathfrak{U},h}^{s,\circ}(N)(-r)) & = & \lim_{\leftarrow j}H^{1}(G_{\mathbf{Q}(\zeta_{mp^{j}})},\mathbf{V}_{\mathfrak{U},h}^{s,\circ}(N)(-r))\\
 & = & H^{1}(G_{\mathbf{Q}(\zeta_{m})},\mathbf{V}_{\mathfrak{U},h}^{s,\circ}(N)(-r)\otimes_{\sco(\mathfrak{U})}\sco(\mathfrak{U})\widehat{\otimes}\mathbf{Z}_{p}[[\Gamma_{m}]])^{\iota}.\end{eqnarray*}
By Proposition \ref{eulerbasic}, the inverse system $\left(_{c,d}\mathfrak{z}_{N,mp^{j}}^{s}(\mathfrak{U},r,a(A))_{h}\right)_{j\geq1}$
defines an element\[
_{c,d}\mathfrak{Z}_{N,m}^{s}(\mathfrak{U},r,a(A))_{h}\in H_{\mathrm{Iw}}^{1}(G_{\mathbf{Q}(\zeta_{m})},\mathbf{V}_{\mathfrak{U},h}^{s,\circ}(N)(-r)).\]

\begin{prop}
The classes $_{c,d}\mathfrak{Z}_{N,m}^{s}(\mathfrak{U},r,a(A))_{h}$
satisfy the compatibility\[
\mathrm{Tw}_{r_{1}-r_{2}}\left(_{c,d}\mathfrak{Z}_{N,m}^{s}(\mathfrak{U},r_{1},a(A))_{h}\right)={}_{c,d}\mathfrak{Z}_{N,m}^{s}(\mathfrak{U},r_{2},a(A))_{h}\]
under the twisting isomorphism\[
\mathrm{Tw}_{r_{1}-r_{2}}:H_{\mathrm{Iw}}^{1}(G_{\mathbf{Q}(\zeta_{m})},\mathbf{V}_{\mathfrak{U},h}^{s,\circ}(N)(-r_{1}))\simeq H_{\mathrm{Iw}}^{1}(G_{\mathbf{Q}(\zeta_{m})},\mathbf{V}_{\mathfrak{U},h}^{s,\circ}(N)(-r_{2})).\]

\end{prop}
Inverting $p$ and passing to the inverse limit over $s$ yields classes\[
_{c,d}\mathfrak{z}_{N,m}(\mathfrak{U},r,a(A))_{h}\in H^{1}\left(G_{\mathbf{Q}(\zeta_{m})},\scv_{\mathfrak{U},h}(N)(-r)\right)\]
and\begin{eqnarray*}
_{c,d}\mathfrak{Z}_{N,m}(\mathfrak{U},r,a(A))_{h} & \in & H_{\mathrm{Iw}}^{1}(G_{\mathbf{Q}(\zeta_{m})},\scv_{\mathfrak{U},h}(N)(-r))\\
 & := & H^{1}(G_{\mathbf{Q}(\zeta_{m})},\scv_{\mathfrak{U},h}(N)(-r)\otimes_{\sco(\mathfrak{U}^{\mathrm{rig}})^{b}}\sco(\mathfrak{U}^{\mathrm{rig}}\times\scx_{m})^{b})^{\iota}\end{eqnarray*}
(here $\scx_{m}=\mathrm{Spf}\mathbf{Z}_{p}[[\Gamma_{m}]]$). 

Now everything glues over an admissible cover of the type constructed
in Lemma \ref{goodcov}, and we get global sections $_{c,d}\mathfrak{z}_{N,m}(r,a(A))$
and $_{c,d}\mathfrak{Z}_{N,m}(r,a(A))$ of the appropriate Galois
cohomology sheaves.

Let $\nu$ be as in the introduction. Set \[
\mu_{c,d}(r,\nu)=(c^{2}-c^{r}\chi_{\scw}(c)^{-1}\nu(c)^{-1}\sigma_{c})(d^{2}-d^{r}\nu(d)\sigma_{d})\in\left(\sco(\scw)\otimes_{\mathbf{Q}_{p}}\mathbf{Q}_{p}(\nu)\right)[\mathrm{Gal}(\mathbf{Q}(\zeta_{m})/\mathbf{Q})],\]
and choose $c,d$ such that $c\equiv1\,\mathrm{mod}\, mp$, $d\equiv1\,\mathrm{mod}\, Nmp$,
and $\nu(c)=\nu(d)=-1$; the existence of $c,d$ with these properties
follows immediately from the Chinese remainder theorem together with
our assumption that $\nu$ has even order. The quantity $\mu_{c,d}$
is then a unit in $\left(\sco(\scw)\otimes_{\mathbf{Q}_{p}}\mathbf{Q}_{p}(\nu)\right)[\mathrm{Gal}(\mathbf{Q}(\zeta_{m})/\mathbf{Q})]$
and\[
\mathfrak{z}_{N,m}(r,\nu)=\mu_{c,d}(r,\nu)^{-1}\sum_{a\in(\mathbf{Z}/A\mathbf{Z})^{\times}}\nu(a)\phantom{}_{c,d}\mathfrak{z}_{m}(r,a(A))\]
is well-defined independently of choosing $c,d$. We have analogously
defined Iwasawa classes $\mathfrak{Z}_{N,m}(r,\nu)$.

\section{Applications to \emph{p}-adic \emph{L}-functions}

\subsection{The regulator map}

In this subsection we prove Theorem \ref{reg}. Let \[
\srr_{\mathbf{Q}_{p}}=\left\{ f\in\mathbf{Q}_{p}[[T,T^{-1}]]\mid f\,\mathrm{convergent}\,\mathrm{on}\,1>|T|>r\,\mathrm{for}\,\mathrm{some}\, r=r_{f}<1\right\} \]
be the Robba ring over $\mathbf{Q}_{p}$, and set $\srr_{L}=\srr_{\mathbf{Q}_{p}}\otimes_{\mathbf{Q}_{p}}L$
and $\srr_{L}^{+}=\srr_{L}\cap L[[T]]$. For background on $(\varphi,\Gamma)$-modules
over $\srr_{L}$, we refer the reader to \cite{BerTri,CC,KPX,Pottharstcyc}.

Notation as in the theorem, fix an affinoid $U\subset\scc$. We define
$\mathrm{Log}|_{U}$ as follows: For any $\mathbf{z}\in H_{\mathrm{Iw}}^{1}(\mathbf{Q}_{p},\scv_{0}(U))$
and $v\in\scd_{\mathrm{crys}}^{\ast}(U)$, $\mathrm{Log}(\mathbf{z})(v)$
is the image of $\mathbf{z}$ under the sequence of maps\begin{multline*}
H_{\mathrm{Iw}}^{1}(\mathbf{Q}_{p},\scv_{0}(U))\hookrightarrow H_{\mathrm{Iw}}^{1}(\mathbf{Q}_{p},\scv_{0}(U))\otimes_{\Lambda(U)}\sco(U\times\scx)\overset{\sim}{\to}\mathbf{D}_{\mathrm{rig}}^{\dagger}(\scv_{0}(U))^{\psi=1}\overset{v^{\ast}}{\longrightarrow}\srr(U)^{+,\psi=(\alpha^{c})^{-1}}\\
\overset{1-\varphi\circ\psi}{\longrightarrow}\srr(U)^{+,\psi=0}\overset{\scm^{-1}}{\longrightarrow}\sco(U)\widehat{\otimes}\mathbf{D}_{\mathrm{an}}(\Gamma)\overset{\sca}{\to}\sco(U\times\scx).\end{multline*}

Let us explain these arrows one-by-one. The first arrow follows from
the identification $\Lambda(U)\cong\sco(U\times\scx)^{b}$ together
with the faithful flatness of $\sco(U\times\scx)$ over $\sco(U\times\scx)^{b}$.
The second arrow follows upon composing the canonical isomorphisms\begin{eqnarray*}
H_{\mathrm{Iw}}^{1}(\mathbf{Q}_{p},\scv_{0}(U))\otimes_{\Lambda(U)}\sco(U\times\scx) & \cong & H_{\mathrm{an.Iw}}^{1}(\mathbf{Q}_{p},\scv_{0}(U))\\
 & \cong & \mathbf{D}_{\mathrm{rig}}^{\dagger}(\scv_{0}(U))^{\psi=1}.\end{eqnarray*}
We refer to Corollary 4.4.11 of  \cite{KPX} for details.

For the third arrow, recall that \begin{eqnarray*}
\scd_{\mathrm{crys}}^{\ast}(U) & = & \mathbf{D}_{\mathrm{rig}}^{\dagger}(\scv_{0}^{\ast}(U))^{\Gamma=1,\varphi=\alpha^{c}}\\
 & \subset & \mathbf{D}_{\mathrm{rig}}^{\dagger}(\scv_{0}^{\ast}(U))\\
 & \cong & \mathrm{Hom}_{\srr(U)}\left(\mathbf{D}_{\mathrm{rig}}^{\dagger}(\scv_{0}(U)),\srr(U)\right),\end{eqnarray*}
so $v$ defines a linear functional $v^{\ast}$ into $\srr(U)$. To
calculate the $\psi$-action on the image of this functional evaluated
on some $d\in\mathbf{D}_{\mathrm{rig}}^{\dagger}(\scv_{0}(U))^{\psi=1}$,
note that the $\varphi$-action on this Hom space is characterized
by the equation $(\varphi\cdot f)(\varphi(d))=\varphi(f(d))$, and
likewise the $\psi$-action is characterized by $(\psi\cdot f)(d)=\psi(f(\varphi(d)))$.
Let $f_{v}=v^{\ast}$, so in particular $\varphi(f_{v}(d))=\alpha^{c}f_{v}(\varphi(d))$
and $f_{v}(d)=\alpha^{c}\psi(f_{v}(\varphi(d)))$. By the definition
of an étale $\varphi$-module over $\srr(U)$ we may (locally on $U$)
write $d=\sum_{i\in I}a_{i}\varphi(d_{i})$ for some $d_{i}\in\mathbf{D}_{\mathrm{rig}}^{\dagger}(\scv_{0}(U))$
and $a_{i}\in\srr(U)$, and since $\psi(d)=d$ we have $\sum_{i\in I}a_{i}\varphi(d_{i})=\sum_{i\in I}\psi(a_{i})d_{i}$.
Now we calculate\begin{eqnarray*}
\psi\left(f_{v}(d)\right) & = & \psi\left(f_{v}\left(\sum a_{i}\varphi(d_{i})\right)\right)\\
 & = & \psi\left(\sum a_{i}f_{v}(\varphi(d_{i}))\right)\\
 & = & (\alpha^{c})^{-1}\psi\left(\sum a_{i}\varphi(f_{v}(d_{i}))\right)\\
 & = & (\alpha^{c})^{-1}\sum\psi(a_{i})f_{v}(d_{i})\\
 & = & (\alpha^{c})^{-1}f_{v}\left(\sum\psi(a_{i})d_{i}\right)\\
 & = & (\alpha^{c})^{-1}f_{v}(d),\end{eqnarray*}
as claimed.

To check that the third arrow really has image in the claimed subspace,
note that for $L/\mathbf{Q}_{p}$ finite and any $a\in\mathcal{O}_{L}$,
there is a natural inclusion $\srr_{L}^{\psi=a^{-1}}\subset\srr_{L}^{+}$
unless $a=1$, in which case $\srr_{L}^{\psi=1}\subset\srr_{L}^{+}\oplus L\cdot\frac{1}{T}\subset\srr_{L}$.
(This is a special case of \cite{ColSPU}, Prop. I.11.) In particular,
given a reduced affinoid $U$ and a function $f\in\sco(U)^{\circ}$
with $f\neq1$, there is an inclusion\[
\left(\srr_{\mathbf{Q}_{p}}\widehat{\otimes}\sco(U)\right)^{\psi=f^{-1}}\subset\srr_{\mathbf{Q}_{p}}^{+}\widehat{\otimes}\sco(U).\]
Now at any crystalline point $x$ of weight $k_{x}\in\mathbf{Z}_{\geq2}$,
the specialization $\alpha_{x}^{c}$ is a \emph{p-}Weil number of
weight $k_{x}-1$, and such points are dense in $\scc$, so the zero
locus of $\alpha^{c}-1$ is nowhere dense in $\scc$.

The fifth arrow is the inverse of the Mellin transform isomorphism\[
\scm:\mathbf{D}_{\mathrm{an}}(\Gamma)\widehat{\otimes}\sco(U)\overset{\sim}{\to}\srr(U)^{+,\psi=0}\]
induced by the isomorphism\begin{eqnarray*}
\mathbf{D}_{\mathrm{an}}(\Gamma) & \overset{\sim}{\to} & \srr_{\mathbf{Q}_{p}}^{+,\psi=0}\\
\mu & \mapsto & \int(1+T)^{\chi_{\mathrm{cyc}}(\gamma)}\mu(\gamma),\end{eqnarray*}
where $\mathbf{D}_{\mathrm{an}}(\Gamma)$ denotes the ring of locally
analytic distributions on $\Gamma$.

The sixth arrow is the Amice transform.

The comparison with the Perrin-Riou regulator map follows from Berger's
alternate construction \cite{BerBK} of the latter (cf. the discussion
following Theorem A in \cite{LLZ}).

It remains to show the claimed growth property at individual points. To
see this, we note that for any point $x \in U( \overline{\mathbf{Q}_p}) $, we have a commutative
square\[
\xymatrix{H_{\mathrm{Iw}}^{1}(\mathbf{Q}_{p},\scv_{0,x})\ar[d]^{\wr}\ar[r] & H_{\mathrm{Iw}}^{1}(\mathbf{Q}_{p},\scv_{0,x})\otimes_{\sco(\scx)^{b}}\sco(\scx)\ar[d]^{\wr}\\
\mathbf{D}^{\dagger}(\scv_{0,x})\ar[r] & \mathbf{D}_{\mathrm{rig}}^{\dagger}(\scv_{0,x})}
\]
of injective maps; here the top and righthand maps are the specialization
at $x$ of the first two maps in the definition of $\mathrm{Log}$,
and the lefthand isomorphism is a well-known theorem of Fontaine (cf. \cite{CC}).
 Going around the lower left half of the square, the claim about growth
now follows from the following result, applied to the representation
$\scv_{0,x}$.
\begin{prop}
Let $V$ be a two-dimensional trianguline representation of $G_{\mathbf{Q}_{p}}$
on an $L$-vector space, and suppose $\mathbf{D}_{\mathrm{rig}}^{\dagger}(V)$
admits a triangulation of the form\[
0\to\srr_{L}(\eta(x_{0})\mu_{\alpha})\to\mathbf{D}_{\mathrm{rig}}^{\dagger}(V)\to\srr_{L}(\mu_{\beta})\to0,\]
where $\mu_{\alpha}$ and $\mu_{\beta}$ are unramified characters
of $\mathbf{Q}_{p}^{\times}$ and $\eta$ is some character of $\mathbf{Z}_{p}^{\times}$.
$ $Let $v_{\alpha},v_{\beta}$ be a basis of $\mathbf{D}_{\mathrm{rig}}^{\dagger}(V)$
realizing this triangulation. Then given any element $f\in\mathbf{D}^{\dagger}(V)\subset\mathbf{D}_{\mathrm{rig}}^{\dagger}(V)$,
the image of $f$ in $\srr_{L}(\mu_{\beta})$ is of the form $v_{\beta}\cdot f_{\beta}$
where $f_{\beta}\in\srr_{L}$ has growth of order $\leq v(\beta^{-1})$.
\end{prop}
\emph{Proof. }This follows from (the proof of) Proposition 3.11 of \cite{ColTri2}.
$\square$

\subsection{Specialization at noble points}

In this section we analyze the specializations of the classes $\mathfrak{z}_{m}(r,\nu)$
and $\mathfrak{Z}_{m}(r,\nu)$ at noble points and prove Theorem \ref{twovarlfn}
and Theorem \ref{critical}. 

Let $f\in S_{k}(\Gamma_{1}(N))$ be a cuspidal newform, $\alpha$
a root of the $p$th Hecke polynomial of $f$, and $\mathbf{f}$ the
associated \emph{p-}stabilization of $f$. We assume that $x_{\mathbf{f}}$
is noble. Set \[
V_{\mathbf{f}}=\left(H_{\mathrm{et}}^{1}(Y_{1}(Np)_{\overline{\mathbf{Q}}},\scl_{k})(2-k)\otimes_{\mathbf{Q}_{p}}\mathbf{Q}_{p}(\mathbf{f})\right)[\ker\phi_{f,\alpha}]\]
and $V_{\mathbf{f}}^{\circ}=$obvious lattice, so $V_{\mathbf{f}}(k)$
is simply the fiber of $\scv$ at $x_{\mathbf{f}}$. Define\begin{eqnarray*}
H_{\mathrm{Iw}}^{1}(V_{\mathbf{f}}^{\circ}) & = & \lim_{\leftarrow j}H^{1}\left(G_{\mathbf{Q}(\zeta_{p^{j}}),Np\infty},V_{\mathbf{f}}^{\circ}\otimes(\mathbf{Z}/p^{j})\right)\end{eqnarray*}
and $H_{\mathrm{Iw}}^{1}(V_{\mathbf{f}})=H_{\mathrm{Iw}}^{1}(V_{\mathbf{f}}^{\circ})[\frac{1}{p}]$.
Recall that $\scx=\mathrm{Spf}(\mathbf{Z}_{p}[[\mathrm{Gal}(\mathbf{Q}(\zeta_{p^{\infty}})/\mathbf{Q})]])^{\mathrm{an}}$,
so $\sco(\scx)^{\circ}=\mathbf{Z}_{p}[[\mathrm{Gal}(\mathbf{Q}(\zeta_{p^{\infty}})/\mathbf{Q})]]$.
\begin{prop}
The module $H_{\mathrm{Iw}}^{1}(V_{f})$ is locally free of rank one
over $\sco(\scx)^{b}$.
\end{prop}
\emph{Proof. }The ring $\sco(\scx)^{b}$ is a one-dimensional, Noetherian,
regular, Jacobson ring, and in particular every prime ideal is maximal
and principal. If $\mathfrak{m}$ is a maximal ideal, then the usual
long exact sequence in cohomology gives a surjection\[
H^{0}(G_{\mathbf{Q},Np\infty},V_{f}\otimes_{\mathbf{Q}_{p}}\sco(\scx)^{b}/\mathfrak{m})\twoheadrightarrow H_{\mathrm{Iw}}^{1}(V_{f})[\mathfrak{m}],\]
and the $H^{0}$ vanishes by the irreducibility of $V_{\mathbf{f}}$.
Therefore $H_{\mathrm{Iw}}^{1}(V_{\mathbf{f}})$ is torsion-free,
so is locally free of some rank. The latter is one by Tate's global
Euler characteristic formula and Kato's results on $H_{\mathrm{Iw}}^{2}$
(cf. Theorem 12.4 of \cite{Kato}). $\square$

Kato's 'optimal Euler system' is an injective $\mathbf{Q}_{p}(f)$-linear
map\begin{eqnarray*}
V_{\mathbf{f}} & \to & H_{\mathrm{Iw}}^{1}(V_{\mathbf{f}}).\\
\gamma & \mapsto & \mathbf{z}_{\gamma,\mathbf{f}}^{(p)}\end{eqnarray*}
with a number of remarkable properties. Let us note that, strictly
speaking, Kato constructs an analogous map $V_{f}\to H_{\mathrm{Iw}}^{1}(V_{f})$
for the unrefined Galois representation $V_{f}\subset H_{\mathrm{et}}^{1}(Y_{1}(N)_{\overline{\mathbf{Q}}},\scl_{k})(2-k)\otimes_{\mathbf{Q}_{p}}\mathbf{Q}_{p}(f)$,
but his construction extends verbatim to nobly refined forms.
\begin{prop}
\label{iwasawaspecialize}$\phantom{}$

\emph{i. }We have\[
\tfrac{1-\nu(-1)\cdot\sigma_{-1}}{2}\mathfrak{Z}_{1}(k,\nu)_{x_{\mathbf{f}}}=\prod_{\ell|A}P_{\ell}(\ell^{-k}\sigma_{\ell}^{-1})_{x_{\mathbf{f}}}\cdot\mathbf{z}_{\delta(\mathbf{f},\nu)}^{(p)},\]
where $\delta(\mathbf{f},\nu)$ is the projection to the $\mathbf{f}$-eigenspace
of $V_{k}(Np)$ of the cohomology class $$\sum_{a\in(\mathbf{Z}/A)^{\times}}\nu(a)\delta_{1,Np}(k,k-1,a(A)).$$

\emph{ii. }Under the hypothesis of Theorem \ref{euler1}.iii, the
class $\delta(\mathbf{f},\nu)$ is nonzero.
\end{prop}
\emph{Proof. }This follows from Proposition \ref{compaA} together
with Lemma 13.11 of Kato. In fact, writing $\delta(\mathbf{f},\nu^{\pm})=b^{\pm}\cdot\gamma^{\pm}$
with $\gamma^{\pm}$ any basis of $V_{\mathbf{f}}^{\pm}$, a calculation
shows that $b^{\pm}$ is a nonzero multiple of $L_{\mathrm{alg}}(k-1,f^{c}\otimes\nu^{-1})$
$\square$

\emph{Proof of Theorem \ref{twovarlfn}. }Notation as in the theorem,
write \[
\mathbf{L}_{\nu^{\pm}}=e^{\pm}\prod_{\ell|\frac{A^{\mp}}{(A^{+},A^{-})}}P_{\ell}(\ell^{-1}\sigma_{\ell}^{-1})\cdot\mathrm{Log}(\mathrm{res}_{p}\mathfrak{Z}_{1}(1,\nu^{\pm})),\]
so $\mathbf{L}_{\nu^{+},\nu^{-}}=\mathbf{L}_{\nu^{+}}+\mathbf{L}_{\nu^{-}}$.
Let \[
\mathrm{log}_{\mathbf{f}}:H_{\mathrm{Iw}}^{1}(\mathbf{Q}_{p},V_{\mathbf{f}}(k-1))\to\sco(\scx)\otimes_{\mathbf{Q}_{p}}\mathbf{D}_{\mathrm{crys}}(V_{\mathbf{f}}(k-1))\]
be the Perrin-Riou regulator map. Kato defines $L_{p,\mathrm{alg}}(\mathbf{f}^{c})$
as $\eta^{\ast}\log_{\mathbf{f}}(\mathbf{z}_{\gamma,\mathbf{f}}^{(p)}(k-1))$
for a certain $\eta\in\mathbf{D}_{\mathrm{crys}}(V_{\mathbf{f}}(k-1))^{\vee}\simeq\mathbf{D}_{\mathrm{crys}}(V_{\mathbf{f}^{c}})$
satisfying $\varphi(\eta)=\alpha^{c}\eta$. By Proposition \ref{iwasawaspecialize}
and a simple twisting argument, we see that\begin{eqnarray*}
\prod_{\ell|A^{\pm}}P_{\ell}(\ell^{-1}\sigma_{\ell}^{-1})_{x_{\mathbf{f}}}\cdot\mathbf{z}_{\delta(\mathbf{f},\nu^{\pm})}^{(p)}(k-1) & = & \tfrac{1+\nu(-1)(-1)^{k}\cdot\sigma_{-1}}{2}\mathfrak{Z}_{1}(1,\nu^{\pm})_{x_{\mathbf{f}}}\\
 & \in & H_{\mathrm{Iw}}^{1}(V_{\mathbf{f}}(k-1)).\end{eqnarray*}
Therefore by Theorem \ref{reg}, we have\begin{eqnarray*}
\mathbf{L}_{\nu^{\pm}}(v)_{x_{\mathbf{f}}} & = & e^{\pm}\prod_{\ell|\frac{A^{\mp}}{(A^{+},A^{-})}}P_{\ell}(\ell^{-1}\sigma_{\ell}^{-1})\cdot v_{x_{\mathbf{f}}}^{\ast}\log_{\mathbf{f}}\left(\prod_{\ell|A^{\pm}}P_{\ell}(\ell^{-1}\sigma_{\ell}^{-1})_{x_{\mathbf{f}}}\cdot\mathbf{z}_{\delta(\mathbf{f},\nu^{\pm})}^{(p)}(k-1)\right)\\
 & = & e^{\pm}\prod_{\ell|A}P_{\ell}(\ell^{-1}\sigma_{\ell}^{-1})_{x_{\mathbf{f}}}\cdot v_{x_{\mathbf{f}}}^{\ast}\log_{\mathbf{f}}\left(\mathbf{z}_{\delta(\mathbf{f},\nu^{\pm})}^{(p)}(k-1)\right)\\
 & = & e^{\pm}\prod_{\ell|A}P_{\ell}(\ell^{-1}\sigma_{\ell}^{-1})_{x_{\mathbf{f}}}\cdot b^{\pm}v_{x_{\mathbf{f}}}^{\ast}\log_{\mathbf{f}}\left(\mathbf{z}_{\gamma^{\pm}}^{(p)}(k-1)\right)\\
 & = & C_{\nu^{\pm}}(v_{x_{\mathbf{f}}})e^{\pm}\prod_{\ell|A}P_{\ell}(\ell^{-1}\sigma_{\ell}^{-1})_{x_{\mathbf{f}}}\cdot L_{p,\mathrm{alg}}(\mathbf{f}^{c})\end{eqnarray*}
as claimed. $\square$

\emph{Proof of Theorem \ref{critical}. }Fix a small connected affinoid
$U\subset\scc$ containing $x_{\mathbf{f}}$ and some nowhere vanishing
section $v\in\scd_{\mathrm{crys}}^{\ast}(U)$, and choose $\gamma\in\scv(U)$
with $\gamma_{x_{\mathbf{f}}}^{\pm}\neq0$. By our previous analysis
together with a deep theorem of Rohrlich \cite{Rohr}, we may choose
$\nu^{\pm}$ such that $ $$\mathbf{L}_{\nu^{\pm}}$ is not identically
zero at $x_{\mathbf{f}}$. Let $\scx^{\pm}$ denote either portion
of $\scx$, and consider the ratio\[
R^{\pm}=\frac{c^{\ast}\mathbf{L}_{\nu^{\pm}}(v)}{\mathbf{L}_{p,\mathrm{an}}^{\pm}\cdot c^{\ast}\prod_{\ell|A^{+}A^{-}}P_{\ell}(\ell^{-1}\sigma_{\ell}^{-1})}\in\mathrm{Frac}(\sco(U\times\scx^{\pm})).\]
By the same theorem of Rohrlich, $L_{p,\mathrm{an}}(\mathbf{f})|_{\scx_{i}}\neq0$
on any connected component $\scx_{i}$ of $\scx$, so the denominator
is not a zero-divisor and this ratio is well-defined. Shrinking $U$
if necessary, we may assume the denominator of $R^{\pm}$ is a non-zero-divisor
after specialization at any noble point $y\in U(\overline{\mathbf{Q}_{p}})$.
Choosing any such point $y$ of noncritical slope, we have \begin{eqnarray*}
R_{y}^{\pm} & = & \frac{\mathbf{L}_{\nu^{\pm}}(v)_{y^{c}}}{L_{p,\mathrm{an}}^{\pm}(\mathbf{f}_{y})\cdot\prod_{\ell|A^{+}A^{-}}P_{\ell}(\ell^{-1}\sigma_{\ell}^{-1})_{y^{c}}}\\
 & = & \frac{C_{\nu^{\pm}}(v)_{y^{c}}\cdot L_{p,\mathrm{alg}}^{\pm}(\mathbf{f}_{y})\cdot\prod_{\ell|A^{+}A^{-}}P_{\ell}(\ell^{-1}\sigma_{\ell}^{-1})_{y^{c}}}{L_{p,\mathrm{an}}^{\pm}(\mathbf{f}_{y})\cdot\prod_{\ell|A^{+}A^{-}}P_{\ell}(\ell^{-1}\sigma_{\ell}^{-1})_{y^{c}}}\\
 & = & C_{y}^{\pm}\in E_{y}\end{eqnarray*}
for some constant $C_{y}$ by Visik's theorem. By the Zariski-density
of such points, we deduce that \[
R^{\pm}\in\mathrm{Frac}(\sco(U))\subset\mathrm{Frac}(\sco(U\times\scx^{\pm})).\]
Since $L_{p,\mathrm{an}}(\mathbf{f})$ is not a zero-divisor, the
polar divisor of $R^{\pm}$ is necessarily disjoint from $\{x_{\mathbf{f}}\}\times\scx^{\pm}$,
so we may specialize $R$ at $x_{\mathbf{f}}$ and conclude that $R_{x_{\mathbf{f}}}^{\pm}$
is a constant, and is furthermore nonzero, since $C_{\nu^{\pm}}(v)_{x_{\mathbf{f}}^{c}}$
is a nonzero multiple of $L(k-1,f\otimes\nu^{\pm})$. Since furthermore
$\prod_{\ell|A^{+}A^{-}}P_{\ell}(\ell^{-1}\sigma_{\ell}^{-1})_{y^{c}}$
is not a zero-divisor, we deduce that $L_{p,\mathrm{an}}^{\pm}(\mathbf{f})=C_{x_{\mathbf{f}}}^{\pm}L_{p,\mathrm{alg}}(\mathbf{f})$,
and evaluating either side at any interpolatory point $x^{j}\eta(x)\in\scx^{\pm}$
for which $ $$L(j+1,f\otimes\eta^{-1})\neq0$ implies $C_{x_{\mathbf{f}}}^{\pm}=1$.
$\square$

\bibliographystyle{amsalpha}

\begin{thebibliography}{MSD74}

\bibitem[AIS13]{AIS}
Fabrizio Andreatta, Adrian Iovita, and Glenn Stevens, \emph{Overconvergent
  {E}ichler-{S}himura isomorphisms}, to appear in J. Inst. Math. Jussieu.

\bibitem[Bel12]{BellCrit}
Jo{\"e}l Bella{\"{\i}}che, \emph{Critical {$p$}-adic {$L$}-functions}, Invent.
  Math. \textbf{189} (2012), no.~1, 1--60. \MR{2929082}

\bibitem[Ber03]{BerBK}
Laurent Berger, \emph{Bloch and {K}ato's exponential map: three explicit
  formulas}, Doc. Math. (2003), no.~Extra Vol., 99--129 (electronic), Kazuya
  Kato's fiftieth birthday. \MR{2046596 (2005f:11268)}

\bibitem[Ber11]{BerTri}
\bysame, \emph{Trianguline representations}, Bull. Lond. Math. Soc. \textbf{43}
  (2011), no.~4, 619--635. \MR{2820149 (2012h:11079)}

\bibitem[BGR84]{BGR}
S.~Bosch, U.~G{\"u}ntzer, and R.~Remmert, \emph{Non-{A}rchimedean analysis},
  Grundlehren der Mathematischen Wissenschaften [Fundamental Principles of
  Mathematical Sciences], vol. 261, Springer-Verlag, Berlin, 1984, A systematic
  approach to rigid analytic geometry. \MR{746961 (86b:32031)}

\bibitem[BK90]{BK}
Spencer Bloch and Kazuya Kato, \emph{{$L$}-functions and {T}amagawa numbers of
  motives}, The {G}rothendieck {F}estschrift, {V}ol.\ {I}, Progr. Math.,
  vol.~86, Birkh\"auser Boston, Boston, MA, 1990, pp.~333--400. \MR{1086888
  (92g:11063)}

\bibitem[Buz07]{BuEigen}
Kevin Buzzard, \emph{Eigenvarieties}, {$L$}-functions and {G}alois
  representations, London Math. Soc. Lecture Note Ser., vol. 320, Cambridge
  Univ. Press, Cambridge, 2007, pp.~59--120. \MR{2392353 (2010g:11076)}

\bibitem[CC99]{CC}
Fr{\'e}d{\'e}ric Cherbonnier and Pierre Colmez, \emph{Th\'eorie d'{I}wasawa des
  repr\'esentations {$p$}-adiques d'un corps local}, J. Amer. Math. Soc.
  \textbf{12} (1999), no.~1, 241--268. \MR{1626273 (99g:11141)}

\bibitem[CE98]{CE}
Robert~F. Coleman and Bas Edixhoven, \emph{On the semi-simplicity of the {$U\sb
  p$}-operator on modular forms}, Math. Ann. \textbf{310} (1998), no.~1,
  119--127. \MR{1600034 (99b:11043)}

\bibitem[CM98]{CMeigencurve}
R.~Coleman and B.~Mazur, \emph{The eigencurve}, Galois representations in
  arithmetic algebraic geometry ({D}urham, 1996), London Math. Soc. Lecture
  Note Ser., vol. 254, Cambridge Univ. Press, Cambridge, 1998, pp.~1--113.
  \MR{1696469 (2000m:11039)}

\bibitem[Col97]{Colemanbanachfamilies}
Robert~F. Coleman, \emph{{$p$}-adic {B}anach spaces and families of modular
  forms}, Invent. Math. \textbf{127} (1997), no.~3, 417--479. \MR{1431135
  (98b:11047)}

\bibitem[Col04]{ColKato}
Pierre Colmez, \emph{La conjecture de {B}irch et {S}winnerton-{D}yer
  {$p$}-adique}, Ast\'erisque (2004), no.~294, ix, 251--319. \MR{2111647
  (2005i:11080)}

\bibitem[Col08]{ColTri2}
\bysame, \emph{Repr\'esentations triangulines de dimension 2}, Ast\'erisque
  (2008), no.~319, 213--258, Repr{\'e}sentations $p$-adiques de groupes
  $p$-adiques. I. Repr{\'e}sentations galoisiennes et $(\phi,\Gamma)$-modules.
  \MR{2493219 (2010f:11173)}

\bibitem[Col10]{ColSPU}
\bysame, \emph{La s\'erie principale unitaire de {${\rm GL}\sb 2(\bold Q\sb
  p)$}}, Ast\'erisque (2010), no.~330, 213--262. \MR{2642407 (2011g:22026)}

\bibitem[Han14]{HanThesis}
David Hansen, \emph{Universal eigenvarieties, trianguline {G}alois
  representations, and {$p$}-adic {L}anglands functoriality}, J. reine angew.
  Math, to appear.

\bibitem[Han15]{Htwovar}
\bysame, \emph{Iwasawa theory of overconvergent modular forms, II: A
  two-variable main conjecture}, in preparation.

\bibitem[Kat04]{Kato}
Kazuya Kato, \emph{{$p$}-adic {H}odge theory and values of zeta functions of
  modular forms}, Ast\'erisque (2004), no.~295, ix, 117--290, Cohomologies
  $p$-adiques et applications arithm{\'e}tiques. III. \MR{2104361
  (2006b:11051)}

\bibitem[Kol88]{Kofinite}
V.~A. Kolyvagin, \emph{Finiteness of {$E({\bf Q})$} and {SH{$(E,{\bf Q})$}} for
  a subclass of {W}eil curves}, Izv. Akad. Nauk SSSR Ser. Mat. \textbf{52}
  (1988), no.~3, 522--540, 670--671. \MR{954295 (89m:11056)}

\bibitem[Kol90]{KoEuler}
\bysame, \emph{Euler systems}, The {G}rothendieck {F}estschrift, {V}ol.\ {II},
  Progr. Math., vol.~87, Birkh\"auser Boston, Boston, MA, 1990, pp.~435--483.
  \MR{1106906 (92g:11109)}

\bibitem[KPX14]{KPX}
Kiran~S. Kedlaya, Jonathan Pottharst, and Liang Xiao, \emph{Cohomology of
  arithmetic families of {$(\varphi,\Gamma)$}-modules}, J. Amer. Math. Soc.
  \textbf{27} (2014), no.~4, 1043--1115. \MR{3230818}

\bibitem[Liu14]{Liu}
Ruochuan Liu, \emph{Triangulation of refined families}, preprint.

\bibitem[LLZ11]{LLZ}
Antonio Lei, David Loeffler, and Sarah Zerbes, \emph{Coleman maps and the
  {$p$}-adic regulator}, Algebra \& Number Theory \textbf{5} (2011), no.~8,
  1095--1131.

\bibitem[LLZ13]{LLZCrit}
Antonio Lei, David Loeffler, and Sarah~Livia Zerbes, \emph{Critical slope
  {$p$}-adic {$L$}-functions of {CM} modular forms}, Israel J. Math.
  \textbf{198} (2013), no.~1, 261--282. \MR{3096640}

\bibitem[MSD74]{MazurSW}
B.~Mazur and P.~Swinnerton-Dyer, \emph{Arithmetic of {W}eil curves}, Invent.
  Math. \textbf{25} (1974), 1--61. \MR{0354674 (50 \#7152)}

\bibitem[Och03]{Och1}
Tadashi Ochiai, \emph{A generalization of the {C}oleman map for {H}ida
  deformations}, Amer. J. Math. \textbf{125} (2003), no.~4, 849--892.
  \MR{1993743 (2004j:11050)}

\bibitem[Och05]{Och2}
\bysame, \emph{Euler system for {G}alois deformations}, Ann. Inst. Fourier
  (Grenoble) \textbf{55} (2005), no.~1, 113--146. \MR{2141691 (2006a:11068)}

\bibitem[Och06]{Och3}
\bysame, \emph{On the two-variable {I}wasawa main conjecture}, Compos. Math.
  \textbf{142} (2006), no.~5, 1157--1200. \MR{2264660 (2007i:11146)}

\bibitem[Pot12]{Pottharstcyc}
Jonathan Pottharst, \emph{Cyclotomic {I}wasawa theory of motives}, preprint.

\bibitem[Pot13]{PottharstAFFSSG}
\bysame, \emph{Analytic families of finite-slope {S}elmer groups}, Algebra
  Number Theory \textbf{7} (2013), no.~7, 1571--1612. \MR{3117501}

\bibitem[PS11]{PSovercvgt}
Robert Pollack and Glenn Stevens, \emph{Overconvergent modular symbols and
  {$p$}-adic {$L$}-functions}, Ann. Sci. \'Ec. Norm. Sup\'er. (4) \textbf{44}
  (2011), no.~1, 1--42. \MR{2760194 (2012m:11074)}

\bibitem[PS13]{PScrit}
\bysame, \emph{Critical slope {$p$}-adic {$L$}-functions}, J. Lond. Math. Soc.
  (2) \textbf{87} (2013), no.~2, 428--452. \MR{3046279}

\bibitem[Roh84]{Rohr}
David~E. Rohrlich, \emph{On {$L$}-functions of elliptic curves and cyclotomic
  towers}, Invent. Math. \textbf{75} (1984), no.~3, 409--423. \MR{735333
  (86g:11038b)}

\bibitem[Rub91]{RuImquad}
Karl Rubin, \emph{The ``main conjectures'' of {I}wasawa theory for imaginary
  quadratic fields}, Invent. Math. \textbf{103} (1991), no.~1, 25--68.
  \MR{1079839 (92f:11151)}

\bibitem[Rub00]{Rubook}
\bysame, \emph{Euler systems}, Annals of Mathematics Studies, vol. 147,
  Princeton University Press, Princeton, NJ, 2000, Hermann Weyl Lectures. The
  Institute for Advanced Study. \MR{1749177 (2001g:11170)}

\bibitem[Ste94]{Strigidsymbs}
Glenn Stevens, \emph{Rigid analytic modular symbols}, unpublished, 1994.

\bibitem[Vi{\v{s}}76]{Visik}
M.~M. Vi{\v{s}}ik, \emph{Nonarchimedean measures associated with {D}irichlet
  series}, Mat. Sb. (N.S.) \textbf{99(141)} (1976), no.~2, 248--260, 296.
  \MR{0412114 (54 \#243)}

\bibitem[Wan12]{W}
Shanwen Wang, \emph{Le syst{\'e}me d'{E}uler de {K}ato en famille, {I}},
  preprint.

\end{thebibliography}
\def\cprime{$'$}
\providecommand{\bysame}{\leavevmode\hbox to3em{\hrulefill}\thinspace}
\providecommand{\MR}{\relax\ifhmode\unskip\space\fi MR }
\providecommand{\MRhref}[2]{%
  \href{http://www.ams.org/mathscinet-getitem?mr=#1}{#2}
}
\providecommand{\href}[2]{#2}

\end{document}